\documentclass[11pt,a4paper,reqno]{amsart}
\usepackage[T1]{fontenc}
\usepackage{lmodern,amsmath,amsfonts,amssymb,amsthm,mathrsfs,mathtools,esint}
\usepackage{microtype,enumitem,comment,needspace}
\usepackage[dvipsnames]{xcolor}
\usepackage[colorlinks=true,linkcolor=blue,citecolor=blue,urlcolor=blue,filecolor=blue,pdfborder={0 0 0},pagebackref=true]{hyperref}
\numberwithin{equation}{section}
\newtheorem{theorem}{Theorem}[section]
\newtheorem{lemma}[theorem]{Lemma}
\newtheorem{proposition}[theorem]{Proposition}
\newtheorem{corollary}[theorem]{Corollary}
\newtheorem{problem}{Problem}
\theoremstyle{definition}
\newtheorem{remark}[theorem]{Remark}

\newcommand{\R}{\mathbb R}
\newcommand{\D}{\mathcal D}
\newcommand{\Spar}{\mathcal S}
\newcommand{\one}{\mathbf 1}
\newcommand{\av}[1]{\langle #1\rangle}
\newcommand{\op}{\mathrm{op}}
\DeclareRobustCommand{\Probref}[1]{\mbox{\hyperref[#1]{Problem~\ref*{#1}}}}
\DeclareRobustCommand{\Tref}[1]{\mbox{\hyperref[#1]{Theorem~\ref*{#1}}}}
\DeclareRobustCommand{\Lref}[1]{\mbox{\hyperref[#1]{Lemma~\ref*{#1}}}}
\DeclareRobustCommand{\Cref}[1]{\mbox{\hyperref[#1]{Corollary~\ref*{#1}}}}
\DeclareRobustCommand{\Sref}[1]{\mbox{\hyperref[#1]{Section~\ref*{#1}}}}
\renewcommand{\eqref}[1]{\hyperref[#1]{\textup{(\ref*{#1})}}}
\hypersetup{pdftitle={Sharp Fefferman--Stein inequalities for matrix reducing maximal operators},pdfauthor={Xi Cen and Pengfei Lu}}
\renewcommand*{\backref}[1]{}

\renewcommand*{\backrefalt}[4]{\ifnum#1>0\space{\footnotesize\color{blue}#2}\fi}
\title[Sharp bounds for reducing maximal operators]{\bfseries Sharp Fefferman--Stein inequalities for matrix reducing maximal operators}
\author[X. Cen]{Xi Cen}
\address{Xi Cen\\
School of Science\\
China University of Mining and Technology-Beijing\\
Beijing 100083, People's Republic of China}
\email{xicenmath@gmail.com}
\author[P. Lu]{Pengfei Lu}
\address{Pengfei Lu\\
Department of Mathematics\\
Northeastern University\\
Shenyang 110004\\
People's Republic of China}
\email{pflumath@gmail.com}
\date{}
\newcommand{\bM}{\overrightarrow{\widetilde{\mathcal M}}}
\newcommand{\bP}{\overrightarrow{\mathcal M}}
\DeclareRobustCommand{\Pref}[1]{\mbox{\hyperref[#1]{Proposition~\ref*{#1}}}}
\DeclareRobustCommand{\Remref}[1]{\mbox{\hyperref[#1]{Remark~\ref*{#1}}}}
\DeclareRobustCommand{\Subref}[1]{\mbox{\hyperref[#1]{Subsection~\ref*{#1}}}}
\begin{document}

\begin{abstract}
We establish sharp Fefferman--Stein inequalities for the matrix reducing maximal operator, answering the question posed by Kakaroumpas and Soler i Gibert~\cite[Math. Z., 2026, before Theorem~6]{KS} concerning the optimal power of the matrix $A_p$ characteristic. For $1<p<\infty$ and $1<q\le\infty$, the norm of the coordinatewise extension from $L^p(\ell^q(\mathbb C^d))$ to $L^p(\ell^q)$ is bounded by $C_{n,d,p}q'[W]_{A_p}^{1/(p-1)}$, with the weight incorporated into the operator and $q'$ the conjugate exponent. The exponent $1/(p-1)$ is optimal. More generally, with the directional Fujii--Wilson characteristic $[W^{-1/(p-1)}]_{A_{p',\infty}^{\mathrm{sc}}}=\sup_{|e|=1}[|W^{-1/p}e|^{p'}]_{A_\infty}$, where $p'=p/(p-1)$ and $e\in\mathbb C^d$ is constant, we prove that the mixed bound by $[W]_{A_p}^{a}[W^{-1/(p-1)}]_{A_{p',\infty}^{\mathrm{sc}}}^b$, $a,b\ge0$, holds uniformly if and only if $a\ge1/p$ and $(p-1)a+b\ge1$. This characterizes all admissible nonnegative mixed powers. The proof uses a common sparse family that controls both the change of reducing matrices and the separate maximal averages of the sequence. On every fixed dyadic grid, our estimate gives a sharp quantitative affirmative answer to Question~7.2 of Domelevo, Kakaroumpas, Petermichl and Soler i Gibert~\cite[J. Math. Anal. Appl., 2024]{DKPS}. For both the full and the dyadic operators, we determine the optimal finite-sequence bound at the identity weight, of order $\min\{q',1+\log N\}$ for $1\le q\le\infty$, with $q'=\infty$ at $q=1$ and $q'=1$ at $q=\infty$.
\end{abstract}

\maketitle

\section{Introduction}\label{sec:intro}
\subsection{Background and motivations}\label{subsec:background}
Fefferman and Stein~\cite[Theorem~1]{FS} extended the classical maximal inequality from a single function to sequences of functions. Let $\mathbb N=\{1,2,\ldots\}$, $n\in\mathbb N$ and let $\mathcal Q$ be the family of bounded open axis-parallel cubes in $\R^n$. For a locally integrable scalar function $h$, the uncentered Hardy--Littlewood maximal operator is
\begin{equation}\label{eq:scalarM}
 Mh(x)=\sup_{\substack{Q\in\mathcal Q\\x\in Q}}\fint_Q|h(y)|\,dy,
\end{equation}
where $\fint_Q h=|Q|^{-1}\int_Qh$ is also denoted by $\av h_Q$. For a scalar sequence $a=(a_j)_j$ and $1<q<\infty$, write $\|a\|_{\ell^q}=(\sum_j|a_j|^q)^{1/q}$. The notation $X\lesssim_\Theta Y$ means $X\le C_\Theta Y$ with a constant depending only on $\Theta$. We write $X\eqsim_\Theta Y$ if both $X\lesssim_\Theta Y$ and $Y\lesssim_\Theta X$ hold. When the parameters are fixed, we may omit the subscript. For a sequence of functions $\mathbf h=(h_j)_j$ and $1<p,q<\infty$, their inequality reads
\[
 \|(Mh_j)_j\|_{L^p(\ell^q)}\lesssim_{n,p,q}\|(h_j)_j\|_{L^p(\ell^q)},\qquad
 \|\mathbf h\|_{L^p(\ell^q)}=\left(\int_{\R^n}\|\mathbf h(x)\|_{\ell^q}^{p}\,dx\right)^{1/p}.
\]

For $1<p<\infty$, Muckenhoupt~\cite{Muck} characterized the locally integrable weights $w>0$ almost everywhere for which $M$ is bounded on $L^p(w)$ by the condition
\begin{equation}\label{eq:scalarAp}
 [w]_{A_p}=\sup_{Q\in\mathcal Q}\left(\fint_Qw\right)\left(\fint_Qw^{-1/(p-1)}\right)^{p-1}<\infty.
\end{equation}
Here $\|f\|_{L^p(w)}=\|w^{1/p}f\|_{L^p}$. The sharp power of $[w]_{A_p}$ in the single-function estimate is $1/(p-1)$ by Buckley~\cite{Buckley}. For sequences the power can be larger. In fact, Cruz-Uribe, Martell and P\'erez~\cite[Theorem~1.12]{CMP} established
\begin{equation}\label{eq:scalar-fs-intro}
 \|(w^{1/p}M(w^{-1/p}h_j))_j\|_{L^p(\ell^q)}\lesssim_{n,p,q}[w]_{A_p}^{\max\{1/q,1/(p-1)\}}\|\mathbf h\|_{L^p(\ell^q)},\qquad 1<p,q<\infty.
\end{equation}
The exponent is optimal, and it equals the single-function exponent precisely when $p\le q+1$. When $p>q+1$, the term $1/q$ records the additional cost of taking separate maximal averages of the sequence. This distinction motivates the question of whether a matrix analogue of the maximal operator necessarily exhibits such a sequence-dependent weight exponent.

A pure $A_p$ power measures the worst possible weight dependence. A more precise bound separates the cost of a single average from the concentration of the weight that persists across scales. For a measurable set $E$ write $\one_E$ for its indicator. For a positive scalar weight $v$, define
\begin{equation}\label{eq:FW}
 [v]_{A_\infty}=\sup_{Q\in\mathcal Q}\frac{1}{v(Q)}\int_Q M(v\one_Q)(x)\,dx,\qquad v(Q)=\int_Qv.
\end{equation}
This is the scalar Fujii--Wilson characteristic used by Hyt\"onen and P\'erez~\cite{HP}. Hyt\"onen and Li~\cite{HyLi} obtained mixed weak- and strong-type estimates for square functions and related operators. Cejas, Li, P\'erez and Rivera-R\'ios~\cite{CLPR} derived mixed weighted estimates for the scalar vector-valued maximal function by sparse domination of the length of the whole sequence. Quantitative extrapolation provides another approach to sequence estimates. For example, Cao, Liu, Si and Yabuta~\cite[Theorem~3.7]{CLSY} obtained weighted square-summable estimates for rough maximal operators.

The preceding estimates concern scalar weights. For functions with values in $\mathbb C^d$, the weight acts on the direction as well as the length of a vector. Volberg~\cite{Volberg} developed the matrix $A_p$ theory, providing the setting in which these two effects must be controlled simultaneously. Fix $d\in\mathbb N$, $1<p<\infty$ and $p'=p/(p-1)$. A matrix weight $W$ on $\R^n$ is a locally integrable $d\times d$ matrix function, positive definite Hermitian almost everywhere. We write $|v|$ for the Euclidean norm on $\mathbb C^d$ and $\|B\|_{\op}=\sup_{|v|=1}|Bv|$. Matrix powers are defined by spectral calculus. In the convention of Kakaroumpas and Soler i Gibert~\cite[(2.3)]{KS}, the characteristic is
\begin{equation}\label{eq:A}
 [W]_{A_p}=\sup_{Q\in\mathcal Q}\fint_Q\left(\fint_Q\|W(x)^{1/p}W(y)^{-1/p}\|_{\op}^{p'}\,dy\right)^{p/p'}dx.
\end{equation}
We say that $W\in A_p$ if this quantity is finite. For $W=wI_d$, where $I_d$ is the identity matrix, \eqref{eq:A} is exactly \eqref{eq:scalarAp}.

To retain the directional information in a mixed estimate, let $Z$ be a matrix weight and $1<s<\infty$. Its scalar directional Fujii--Wilson characteristic at exponent $s$ is
\[
 [Z]_{A_{s,\infty}^{\mathrm{sc}}}=\sup_{\substack{e\in\mathbb C^d\\ |e|=1}}[|Z^{1/s}e|^s]_{A_\infty}.
\]
Here the unit vector $e$ is constant, and the superscript $\mathrm{sc}$ refers to the scalar weights in the supremum. In particular, for the dual matrix weight $W^{-1/(p-1)}$ at exponent $p'$, the directional weights are $|W^{-1/p}e|^{p'}$. The finiteness of the dual characteristic for $W\in A_p$ follows from \Lref{lem:kernel}.

To measure a vector at $x$ after averaging over a cube, Christ and Goldberg~\cite{CG} introduced a matrix maximal operator for $p=2$. Goldberg~\cite{Goldberg} developed the general $p$ theory. In the normalization used throughout this paper, the corresponding pointwise operator is defined for $W\in A_p$ and $h\in L^p(\R^n,\mathbb C^d)$ by
\begin{equation}\label{eq:pointwise-max}
 \mathcal M_W h(x)=\sup_{\substack{Q\in\mathcal Q\\x\in Q}}\fint_Q|W(x)^{1/p}W(y)^{-1/p}h(y)|\,dy.
\end{equation}
At $p=1$, Lerner, Li, Ombrosi and Rivera-R\'ios~\cite[Theorem~1 and Lemma~2]{LLOR} showed the sharp quadratic weight dependence in the corresponding normalized weak-type estimate for the pointwise maximal operator, with scalar examples for the lower bound.

In estimates involving averages of matrix weights over cubes, it is useful to represent the weighted length by a single matrix. For every $Q\in\mathcal Q$, choose a positive definite Hermitian reducing matrix $V_Q$ of $W$ at exponent $p$, so that
\begin{equation}\label{eq:reducing}
 \left(\fint_Q|W(y)^{1/p}v|^p\,dy\right)^{1/p}\le |V_Qv|\le c_d\left(\fint_Q|W(y)^{1/p}v|^p\,dy\right)^{1/p},\quad v\in\mathbb C^d.
\end{equation}
The existence follows from Goldberg~\cite[Proposition~1.2]{Goldberg}. We fix $c_d=\sqrt{2d}$ and require all estimates to hold uniformly over families satisfying \eqref{eq:reducing}. Replacing $W(x)^{1/p}$ in \eqref{eq:pointwise-max} by $V_Q$ gives the reducing maximal operator
\begin{equation}\label{eq:max}
 \widetilde{\mathcal M}_W h(x)=\sup_{\substack{Q\in\mathcal Q\\x\in Q}}\fint_Q|V_QW(y)^{-1/p}h(y)|\,dy.
\end{equation}
Both \eqref{eq:pointwise-max} and \eqref{eq:max} act from $L^p(\R^n,\mathbb C^d)$ to nonnegative scalar functions. Their measurability and the local integrability required in the averages follow from \Lref{lem:domain}. For $W=wI_d$ one may take $V_Q=\av w_Q^{1/p}I_d$. Thus
\begin{equation}\label{eq:scalar}
 \begin{aligned}
 \mathcal M_{wI_d}h(x)&=w(x)^{1/p}M(w^{-1/p}|h|)(x),\\
 \widetilde{\mathcal M}_{wI_d}h(x)&=\sup_{Q\ni x}\av w_Q^{1/p}\fint_Q w(y)^{-1/p}|h(y)|\,dy.
 \end{aligned}
\end{equation}
The first operator uses the weight at the point of observation, whereas the second uses a matrix representing its averaged length over the cube. Even for scalar weights, these are different operators. In particular, the optimal scalar exponent in \eqref{eq:scalar-fs-intro} concerns the first one. The second is the operator studied here. Isralowitz, Kwon and Pott~\cite[p.~252]{IKP} obtained the sharp power $1/(p-1)$ for the single-function reducing maximal estimate, as also recalled in \cite[Section~3.1]{KS}.

For $1\le q\le\infty$, let $\ell_d^q=\ell^q(\mathbb C^d)$ denote the space of sequences of vectors with norm $(\sum_j|h_j|^q)^{1/q}$ for $q<\infty$ and $\sup_j|h_j|$ for $q=\infty$. The corresponding scalar space is $\ell^q=\ell^q(\mathbb C)$. For $1<q<\infty$ set $q'=q/(q-1)$, with $q'=\infty$ at $q=1$ and $q'=1$ at $q=\infty$. For $N\in\mathbb N$, write $\ell_{N,d}^q$ and $\ell_N^q$ for the corresponding spaces of sequences of length $N$, with the same norms. For finite or countable coordinatewise measurable sequences $\mathbf h=(h_j)_j$ in $L^p(\R^n,\ell_d^q)$, define
\begin{equation}\label{eq:vector-extension}
 \bM_W\mathbf h=(\widetilde{\mathcal M}_Wh_j)_j,\qquad \bP_W\mathbf h=(\mathcal M_Wh_j)_j.
\end{equation}
We abbreviate the domain $L^p(\R^n,\ell_d^q)$ and the range $L^p(\R^n,\ell^q)$ by $L^p(\ell_d^q)$ and $L^p(\ell^q)$, respectively, with the usual Lebesgue measure norms
\begin{equation}\label{eq:matrixseq}
 \|\mathbf h\|_{L^p(\ell_d^q)}=\left(\int_{\R^n}\|\mathbf h(x)\|_{\ell_d^q}^{p}\,dx\right)^{1/p}.
\end{equation}
We use $\|T\|_{X\to Y}$ for the least constant in $\|Th\|_Y\le C\|h\|_X$, including for sublinear or nonlinear operators, and allow the value $\infty$. The case $q=p$ follows by summing the single-function estimates. The key difficulty is $1<q<p$, where each coordinate can select a different cube and the $L^p$ norm lies outside the $\ell^q$ norm. Thus control of each coordinate alone does not yield the desired uniform sequence bound.

Our notation incorporates the weight into the operators. In the original notation of Kakaroumpas and Soler i Gibert~\cite[Section~3.1]{KS}, $M_W f$ and $\widetilde M_W f$ average $|W(x)^{1/p}f|$ and $|V_Qf|$, respectively. For $h=W^{1/p}f$, we have
\begin{equation}\label{eq:matrixsingle}
 \widetilde{\mathcal M}_Wh=\widetilde M_W f,\qquad \mathcal M_Wh=M_W f.
\end{equation}
The multiplication $f\mapsto W^{1/p}f$ is an isometry from the original weighted domain, whose norm is $\|(W^{1/p}f_j)_j\|_{L^p(\ell_d^q)}$, onto $L^p(\ell_d^q)$. Thus the least constants in the two formulations coincide.

A principal reason to study the square-summable extension in \eqref{eq:vector-extension} comes from product spaces. For the scalar strong maximal operator, Cao, Xue and Yabuta~\cite[Theorem~2.4]{CXY} obtained vector-valued two-weight inequalities for rectangular averages under bump and covering hypotheses. In the matrix-weighted theory of biparameter Journ\'e operators, Domelevo, Kakaroumpas, Petermichl and Soler i Gibert~\cite[Sections~7--8]{DKPS} needed a different maximal estimate: Haar coefficients in one variable form a sequence, and the reducing maximal operator acts in the other variable. Their Proposition~7.1 gives the desired $\ell^2$ estimate on a fixed dyadic grid for $1<p\le2$. In Question~7.2 they ask about $p>2$ under the same grid-dependent weight condition. The single-grid hypothesis is essential to its application: the weight is assumed to satisfy the matrix $A_p$ condition on that grid, without imposing the full-cube condition. In the context of general bases, Cruz-Uribe and \c{S}irin~\cite{CUS} developed matrix extrapolation under a boundedness assumption on the corresponding Christ--Goldberg maximal operator.

For the family of all cubes, Kakaroumpas and Soler i Gibert~\cite[Theorem~1]{KS} subsequently established the bound
\[
 \|\bP_W\|_{L^p(\ell_d^q)\to L^p(\ell^q)}\lesssim_{n,d,p,q}[W]_{A_p}^{\max\{1/(q-1),1/(p-1)\}},\qquad 1<p,q<\infty,
\]
using convex-set extrapolation. To estimate the reducing operator, they first use $\widetilde{\mathcal M}_Wh\lesssim_{d,p}[W]_{A_p}^{1/p}M(\mathcal M_Wh)$, then apply the unweighted Fefferman--Stein inequality. Their Theorem~6 therefore gives
\begin{equation}\label{eq:original-bound}
 \|\bM_W\|_{L^p(\ell_d^q)\to L^p(\ell^q)}\lesssim_{n,d,p,q}[W]_{A_p}^{\frac1p+\max\{\frac1{q-1},\frac1{p-1}\}},\qquad 1<p,q<\infty.
\end{equation}
In passing from the pointwise operator to the reducing operator this argument introduces the additional power $1/p$. The authors explicitly ask whether the exponent can be improved in \cite[p.~21, immediately before Theorem~6]{KS}. The challenge is therefore both qualitative and quantitative: obtain an estimate under the original single-grid hypothesis of \cite[Question~7.2]{DKPS}, and remove the excess weight powers in \eqref{eq:original-bound}. We resolve both issues by treating the reducing operator directly. The mixed characteristic will also allow us to identify exactly which powers are necessary.

\subsection{The original problems and the main results}\label{subsec:results}
The question preceding Theorem~6 in Kakaroumpas and Soler i Gibert~\cite[p.~21]{KS} asks whether the exponent in \eqref{eq:original-bound} can be improved. In the equivalent normalization \eqref{eq:matrixsingle}, determining the best possible improvement leads to the following formulation.
\begin{problem}[{Kakaroumpas--Soler i Gibert~\cite[p.~21, immediately before Theorem~6]{KS}}]\label{prob:KS}
Fix $n,d\ge1$ and $1<p,q<\infty$. Determine the least $\beta$ such that
\[
 \|\bM_W\|_{L^p(\ell_d^q)\to L^p(\ell^q)}\lesssim_{n,d,p,q}[W]_{A_p}^\beta
\]
for every matrix $A_p$ weight, uniformly over all admissible reducing families and sequence lengths.
\end{problem}

The second problem is posed on a single dyadic grid. A dyadic grid $\D$ is a family of half-open cubes partitioning $\R^n$ at each scale $2^{-k}$, $k\in\mathbb Z$, and nested across scales. The standard grid is $\D_0=\{2^{-k}(m+[0,1)^n):k\in\mathbb Z,\ m\in\mathbb Z^n\}$. We write $[W]_{A_p(\D)}$ for \eqref{eq:A} with $Q\in\D$, and use $W\in A_p(\D)$ when this quantity is finite. The operators $M^\D$, $\widetilde{\mathcal M}_W^\D$ and $\bM_W^\D$ are the respective restrictions of \eqref{eq:scalarM}, \eqref{eq:max} and \eqref{eq:vector-extension} to this grid. Define $[v]_{A_\infty(\D)}$ by \eqref{eq:FW} using $Q\in\D$ and $M^\D$. The directional characteristic $[Z]_{A_{t,\infty}^{\mathrm{sc}}(\D)}$ is defined by the same supremum over constant unit vectors, using $[|Z^{1/t}e|^t]_{A_\infty(\D)}$ in place of the full-cube characteristic.
Proposition~7.1 of Domelevo, Kakaroumpas, Petermichl and Soler i Gibert~\cite{DKPS} gives the estimate
\begin{equation}\label{eq:dkps-known-range}
 \|\bM_W^\D\|_{L^p(\ell_d^2)\to L^p(\ell^2)}\lesssim_{n,d,p}[W]_{A_p(\D)}^{1/(p-1)},\qquad 1<p\le2.
\end{equation}
For $1<p<2$, their reverse H\"older argument gives an exponent $r<p$, so the scalar sequence estimate applies at $2/r>1$. At $p=2$ they sum the single-function bounds. For $p>2$ the exponent $r$ can exceed two, as explained in \cite[after Proposition~7.1]{DKPS}. In the equivalent normalization \eqref{eq:matrixsingle}, their Question~7.2 is the following problem.
\begin{problem}[{Domelevo--Kakaroumpas--Petermichl--Soler i Gibert~\cite[Question~7.2]{DKPS}}]\label{prob:DKPS}
For $p>2$ and $W\in A_p(\D)$, is there a finite $C=C(n,d,p,W)$, independent of the sequence, such that
\begin{equation}\label{eq:original-dyadic}
 \|\bM_W^\D\|_{L^p(\ell_d^2)\to L^p(\ell^2)}\le C?
\end{equation}
\end{problem}

We first answer \Probref{prob:KS} by determining the optimal pure power. We then characterize the admissible mixed powers and apply the estimates to \Probref{prob:DKPS}.

\begin{theorem}[Optimal pure matrix $A_p$ exponent]\label{th:main}
Let $n,d\in\mathbb N$, $1<p<\infty$, $1<q\le\infty$ and $W\in A_p$. Then, uniformly over reducing families satisfying \eqref{eq:reducing} and over sequence lengths,
\begin{equation}\label{eq:sharp-pure}
 \|\bM_W\|_{L^p(\ell_d^q)\to L^p(\ell^q)}
 \lesssim_{n,d,p}q'[W]_{A_p}^{1/(p-1)}.
\end{equation}
Moreover, let $t\ge1$. Then
\begin{equation}\label{eq:extremal-norm}
 \sup_{\substack{W\in A_p\\ [W]_{A_p}\le t}}
 \|\bM_W\|_{L^p(\ell_d^q)\to L^p(\ell^q)}
 \eqsim_{n,d,p,q}t^{1/(p-1)}.
\end{equation}
\end{theorem}

\Tref{th:main} answers \Probref{prob:KS} in full. Compared with \cite[Theorem~6]{KS}, the reduction in the exponent is
\[
 \frac1p+\max\left\{0,\frac1{q-1}-\frac1{p-1}\right\},\qquad 1<q<\infty.
\]
In particular, the new exponent is independent of $q$. The equivalence \eqref{eq:extremal-norm} gives more than a necessary lower bound on possible powers: it identifies the growth of the largest operator norm under the weight constraint $[W]_{A_p}\le t$. In particular, for every $0\le\beta<1/(p-1)$ there is no constant independent of $W$ in the bound by $[W]_{A_p}^\beta$. This follows by dividing \eqref{eq:extremal-norm} by $t^\beta$ and letting $t\to\infty$. A single nonzero coordinate suffices for the lower estimate.

The pure estimate depends only on the matrix $A_p$ characteristic. To distinguish the cost of one reducing average from the concentration of the dual weight across scales, we now use the full directional characteristic introduced in \Subref{subsec:background}. For $W=wI_d$, it reduces to the scalar characteristic $[w^{-1/(p-1)}]_{A_\infty}$. By \Lref{lem:kernel},
\[
 1\le [W^{-1/(p-1)}]_{A_{p',\infty}^{\mathrm{sc}}}\lesssim_{n,d,p}[W]_{A_p}^{1/(p-1)}.
\]
This relation permits a trade-off between the two powers. The next theorem identifies exactly which pairs of nonnegative powers are admissible.

\begin{theorem}[Sharp mixed-characteristic region]\label{th:mixed}
Let $n,d\in\mathbb N$, $1<p<\infty$, $1<q\le\infty$ and $W\in A_p$. Then, uniformly over admissible reducing families and sequence lengths,
\begin{equation}\label{eq:main}
 \|\bM_W\|_{L^p(\ell_d^q)\to L^p(\ell^q)}
 \lesssim_{n,d,p}q'[W]_{A_p}^{1/p}[W^{-1/(p-1)}]_{A_{p',\infty}^{\mathrm{sc}}}^{1/p}.
\end{equation}
Fix $n,d,p,q$ and $a,b\ge0$. Then the estimate
\begin{equation}\label{eq:admissible}
 \|\bM_W\|_{L^p(\ell_d^q)\to L^p(\ell^q)}
 \lesssim_{n,d,p,q,a,b} [W]_{A_p}^a[W^{-1/(p-1)}]_{A_{p',\infty}^{\mathrm{sc}}}^b
\end{equation}
holds uniformly for all $W\in A_p$ and reducing families in \eqref{eq:reducing} if and only if
\begin{equation}\label{eq:region}
 a\ge\frac1p,\qquad (p-1)a+b\ge1.
\end{equation}
\end{theorem}

The two inequalities in \eqref{eq:region} arise at opposite endpoints of the scalar power-weight class. For $0<\varepsilon<1/2$, let $w_-(x)=|x_1|^{-1+\varepsilon}$, $w_+(x)=|x_1|^{(p-1)(1-\varepsilon)}$ and $W_\pm=w_\pm I_d$, with arbitrary positive values on $\{x_1=0\}$. The tests in \Subref{sub:tests} give
\begin{align*}
 [W_-]_{A_p}&\eqsim_p\varepsilon^{-1},\qquad [W_-^{-1/(p-1)}]_{A_{p',\infty}^{\mathrm{sc}}}\eqsim_{n,p}1,\\
 \|\bM_{W_-}\|_{L^p(\ell_d^q)\to L^p(\ell^q)}&\eqsim_{n,d,p,q}\varepsilon^{-1/p},\\
 [W_+]_{A_p}&\eqsim_p\varepsilon^{-(p-1)},\qquad [W_+^{-1/(p-1)}]_{A_{p',\infty}^{\mathrm{sc}}}\eqsim_{n,p}\varepsilon^{-1},\\
 \|\bM_{W_+}\|_{L^p(\ell_d^q)\to L^p(\ell^q)}&\eqsim_{n,d,p,q}\varepsilon^{-1}.
\end{align*}
The first family forces $a\ge1/p$, since the dual characteristic stays bounded. The second forces $(p-1)a+b\ge1$, since the two characteristics grow at different rates. Conversely, the mixed upper bound \eqref{eq:main} and the comparison preceding \Tref{th:mixed} imply every bound in \eqref{eq:region}. The two endpoint mechanisms thus determine the entire region. The following formulation makes the sharpness of its lower boundary explicit.

\begin{corollary}[Optimal trade-off between the weight powers]\label{cor:tradeoff}
Suppose that the hypotheses of \Tref{th:mixed} hold and fix $b\ge0$. Then the least admissible exponent of $[W]_{A_p}$ in \eqref{eq:admissible} is
\begin{equation}\label{eq:tradeoff}
 a_{\min}(b)=\max\left\{\frac1p,\frac{1-b}{p-1}\right\}.
\end{equation}
In particular, let $0\le\theta\le1$. Then
\begin{equation}\label{eq:frontier}
 \|\bM_W\|_{L^p(\ell_d^q)\to L^p(\ell^q)}
 \lesssim_{n,d,p}q'[W]_{A_p}^{(1-\theta)/(p-1)+\theta/p}[W^{-1/(p-1)}]_{A_{p',\infty}^{\mathrm{sc}}}^{\theta/p}.
\end{equation}
Every pair of powers in \eqref{eq:frontier} is minimal for the coordinatewise order on the nonnegative quadrant.
\end{corollary}

\begin{remark}[A lower bound for every weight]\label{rem:sharp-law}
The power weights $W_+$ attain the order of the right-hand weight factor in \eqref{eq:frontier}, for all $0\le\theta\le1$. The family $W_-$ isolates the exponent $1/p$. Beyond these examples, the norm of one reducing average in \Lref{lem:localnorm} yields a lower bound for every matrix weight. Together with \Tref{th:mixed}, for $1<q\le\infty$ we have
\[
 \begin{aligned}
 [W]_{A_p}^{1/p}&\lesssim_{d,p}\|\bM_W\|_{L^p(\ell_d^q)\to L^p(\ell^q)},\\
 \|\bM_W\|_{L^p(\ell_d^q)\to L^p(\ell^q)}&\lesssim_{n,d,p}q'[W]_{A_p}^{1/p}[W^{-1/(p-1)}]_{A_{p',\infty}^{\mathrm{sc}}}^{1/p}.
 \end{aligned}
\]
Let $t_0\ge1$ and suppose $[W^{-1/(p-1)}]_{A_{p',\infty}^{\mathrm{sc}}}\le t_0$. Then
\[
 \|\bM_W\|_{L^p(\ell_d^q)\to L^p(\ell^q)}\eqsim_{n,d,p,q,t_0}[W]_{A_p}^{1/p}.
\]
Thus $1/p$ describes the growth for weights with bounded dual characteristic, whereas $1/(p-1)$ in \eqref{eq:extremal-norm} describes the unrestricted worst case. These are different sharpness statements.
\end{remark}

The following corollary answers \Probref{prob:DKPS} and strengthens \eqref{eq:original-dyadic} to a sharp quantitative bound under the original single-grid hypothesis.
\Needspace{8\baselineskip}
\begin{corollary}[A single dyadic grid]\label{cor:dyadic}
Let $n,d\in\mathbb N$, $1<p<\infty$, $1<q\le\infty$, and let $\D$ be a fixed dyadic grid. Suppose $W\in A_p(\D)$. Then, for every reducing family on $\D$ satisfying \eqref{eq:reducing},
\begin{equation}\label{eq:dyadic-main}
 \begin{aligned}
 \|\bM_W^\D\|_{L^p(\ell_d^q)\to L^p(\ell^q)}
 &\lesssim_{n,d,p}q'[W]_{A_p(\D)}^{1/p}[W^{-1/(p-1)}]_{A_{p',\infty}^{\mathrm{sc}}(\D)}^{1/p}\\
 &\lesssim_{n,d,p}q'[W]_{A_p(\D)}^{1/(p-1)}.
 \end{aligned}
\end{equation}
Moreover, let $t\ge1$. Then
\begin{equation}\label{eq:dyadic-extremal}
 \sup_{\substack{W\in A_p(\D)\\{[W]}_{A_p(\D)}\le t}}
 \|\bM_W^\D\|_{L^p(\ell_d^q)\to L^p(\ell^q)}\eqsim_{n,d,p,q}t^{1/(p-1)}.
\end{equation}
All constants are independent of the grid and of the admissible reducing family. The admissible nonnegative mixed powers for each fixed grid are exactly those in \eqref{eq:region}, with both characteristics restricted to $\D$.
\end{corollary}

In particular, \eqref{eq:dyadic-extremal} gives a sharp quantitative answer to \Probref{prob:DKPS} when $q=2$ and $p>2$. The full range of mixed powers in \Cref{cor:dyadic} shows that this answer remains optimal even when the dual directional characteristic is included. The lower bound from \Lref{lem:localnorm} applies on the same grid. In particular, if $[W^{-1/(p-1)}]_{A_{p',\infty}^{\mathrm{sc}}(\D)}\le t_0$ with $t_0\ge1$, then
\[
 \|\bM_W^\D\|_{L^p(\ell_d^q)\to L^p(\ell^q)}\eqsim_{n,d,p,q,t_0}[W]_{A_p(\D)}^{1/p}.
\]

The pointwise Christ--Goldberg maximal operator is the subject of the separate question in \cite[Question~7.3]{DKPS}. Both \Probref{prob:KS} and \Probref{prob:DKPS} concern the reducing-matrix counterpart.

\begin{corollary}[Finite sequences]\label{cor:length}
Let $1<p<\infty$, $N\in\mathbb N$ and $1\le q\le\infty$. Let $\mathcal G$ denote the family of all cubes or a fixed dyadic grid, and suppose $W\in A_p(\mathcal G)$. Write $\bM_{W,N}^{\mathcal G}$ for the restriction of the sequence operator to cubes in $\mathcal G$ and to $N$ entries. Then
\begin{equation}\label{eq:length}
 \begin{aligned}
 \|\bM_{W,N}^{\mathcal G}\|_{L^p(\ell_{N,d}^q)\to L^p(\ell_N^q)}
 &\lesssim_{n,d,p}\min\{q',1+\log N\}\\
 &\quad\cdot[W]_{A_p(\mathcal G)}^{1/p}[W^{-1/(p-1)}]_{A_{p',\infty}^{\mathrm{sc}}(\mathcal G)}^{1/p}.
 \end{aligned}
\end{equation}
If $W=I_d$ and $V_Q=I_d$, then
\begin{equation}\label{eq:length-sharp}
 \|\bM_{I_d,N}^{\mathcal G}\|_{L^p(\ell_{N,d}^q)\to L^p(\ell_N^q)}\eqsim_{n,p}\min\{q',1+\log N\}.
\end{equation}
In particular, the countable-sequence norm at the identity weight is comparable to $q'$ for $1<q\le\infty$ and is infinite for $q=1$, on each grid as well as for all cubes.
\end{corollary}

Taking $W=I_d$ and letting $N$ grow shows that the factor $q'$ has the correct order as $q\downarrow1$ at the identity weight. At $q=1$ the loss is logarithmic in $N$. These sequence bounds use one fixed weight, whereas the lower bound for the weight exponent uses one fixed coordinate.

For the mixed square-maximal function in \cite[Lemma~8.2(1) and Remark~8.3]{DKPS}, \Cref{cor:dyadic} extends the estimate to $p>2$ on the original product grid. We give the operator and its quantitative bound in \Subref{sub:mixed-application}.

The single-function mixed estimate of Kakaroumpas and Nieraeth~\cite{KN} already implies the full-cube bound for $q\ge p$ by a positive-kernel argument (\Subref{subsec:comparison-proof}). For $1<q<p$, the coordinatewise argument would require a reverse H\"older gain that can fail for large weight characteristic. We choose one sparse family for the whole sequence and keep the separate cube suprema within each group of descendants. Scalar sequence and lattice versions of sparse domination appear in Culiuc, Di Plinio and Ou~\cite[Theorem~1.2]{CDO} and H\"anninen and Lorist~\cite{HL}. Here the stopping rule also controls the change from a reducing matrix to those of its descendants.

\subsection{Organization}\label{subsec:organization}
\Sref{sec:mainproof} states the auxiliary estimates and proves the main theorems, their sharpness and the corollaries. \Sref{sec:aux} proves the auxiliary lemmas and the sparse domination proposition.

\Needspace{8\baselineskip}
\subsection{Notation}\label{subsec:notation}
We collect the remaining notation used in the proofs.
\begin{itemize}[label=\textbullet,leftmargin=*,itemsep=4pt,topsep=4pt,parsep=0pt]
\item \textbf{Cubes and averages.} Let $\mathcal G$ denote $\mathcal Q$ or a fixed dyadic grid $\D$. The side length of a cube $Q$ is $\ell(Q)$. A grid cube and its interior give the same integrals, while the grid suprema use half-open cubes. For $1\le r<\infty$, write
\[
 \av h_{r,Q}=\left(\fint_Q|h|^r\right)^{1/r},\qquad \av h_Q=\fint_Q h.
\]
For vector functions, $|h|$ is the Euclidean length and the ordinary average is taken componentwise.
\item \textbf{Norms and exponents.} Write $e_1,\ldots,e_d$ for the standard basis of $\mathbb C^d$, $B^*$ for the matrix adjoint, and $\|h\|_s=\||h|\|_{L^s(\R^n)}$. The notation $\|h\|_{L^s(E)}$ uses the same convention on measurable sets $E$. All logarithms are natural. For any exponent $1<s<\infty$, let $s'=s/(s-1)$.
\item \textbf{Weight characteristics.} For $\mathcal G=\mathcal Q$ we use the full-cube characteristics. For $\mathcal G=\D$ we restrict the defining suprema to that grid. The notations $[W]_{A_p(\mathcal G)}$ and $[Z]_{A_{s,\infty}^{\mathrm{sc}}(\mathcal G)}$ cover both cases. For $\mathcal G=\mathcal Q$ we omit the family from the subscript. In the norm-growth laws, $t\ge1$ is a bound on the weight characteristic, and $t_0\ge1$ is a fixed bound on the dual directional characteristic. Scalar and matrix $A_s(\mathcal G)$ characteristics are taken at exponent $s\in(1,\infty)$. If $W\in A_p$, each restricted characteristic is bounded by its full-cube counterpart.
\item \textbf{Dual weights and reducing matrices.} Set $U=W^{-1/(p-1)}$ for the dual matrix weight, locally integrable by \Lref{lem:domain}. Write $V_{Z,s,Q}$ for a reducing matrix of $Z$ at exponent $s$, and $V'_Q=V_{U,p',Q}$. Two families satisfying \eqref{eq:reducing} give comparable lengths with factor $c_d$. For grid restrictions one may use the matrices on cube interiors or an independent admissible family.
\item \textbf{Single-cube operators.} For $Q\in\mathcal G$, set
\[
 K_Q(y)=\|V_QW(y)^{-1/p}\|_{\op},\quad B_Qh=\one_Q\av{|V_QW^{-1/p}h|}_Q,
 \quad\overrightarrow B_Q\mathbf h=(B_Qh_j)_j.
\]
\item \textbf{Sequence norms.} Finite and countable sequences are understood as in \eqref{eq:matrixseq}. All measurability and truncation assertions are proved in \Lref{lem:domain}. By Tonelli's theorem and a single-entry test,
\begin{equation}\label{eq:diagonal-norm}
 \|\bM_W\|_{L^p(\ell_d^p)\to L^p(\ell^p)}=\|\widetilde{\mathcal M}_W\|_{L^p(\mathbb C^d)\to L^p}.
\end{equation}
All suprema over empty collections of nonnegative quantities are zero.
\item \textbf{Dyadic trees and forests.} A descendant of $Q\in\D$ is a cube $R\in\D$ with $R\subseteq Q$, including $Q$ itself. A strict descendant satisfies $R\subsetneq Q$. The parent is the unique cube of twice the side length containing $Q$. A finite dyadic tree rooted at $Q$ contains $Q$ and all descendants down to a specified generation. A finite dyadic forest $\mathcal T$ is a union of finitely many such trees with pairwise disjoint roots. Terminal cubes form the last generation of their tree and partition each root. The descendants of any $Q\in\mathcal T$ in the forest form a finite tree rooted at $Q$.
\item \textbf{Local maximal operators.} For $h\in L^1(Q)$ define on $Q$
\[
 \begin{aligned}
 M_Q^\D h(x)&=\sup_{\substack{R\in\D\\x\in R\subseteq Q}}\av{|h|}_R,\\
 M_Q^{\mathcal T}h(x)&=\max_{\substack{R\in\mathcal T\\x\in R\subseteq Q}}\av{|h|}_R.
 \end{aligned}
\]
In the second definition $Q\in\mathcal T$. Both operators vanish off $Q$. For a finite forest $\mathcal T$, the notation $\widetilde{\mathcal M}_W^{\mathcal T}$ denotes \eqref{eq:max} restricted to $R\in\mathcal T$, with empty maximum zero. Its sequence extension is defined coordinatewise.
\item \textbf{Sparse families.} A family $\Spar$ is $1/2$-sparse if there exist pairwise disjoint measurable sets $E_Q\subseteq Q$ with $|E_Q|\ge|Q|/2$ for every $Q\in\Spar$.
\end{itemize}

\section{Proofs of the main results}\label{sec:mainproof}
\subsection{The auxiliary lemmas}\label{subsec:auxstatements}
The kernel estimate below supplies an averaging exponent $r<p$. The sparse bound in \Pref{prop:freeze} applies that exponent to the length of the sequence.

\begin{lemma}[Well-definedness]\label{lem:domain}
Let $W\in A_p(\mathcal G)$, where $\mathcal G$ is the family of all cubes or a fixed dyadic grid. Then $W^{-1/p}$ has operator norm in $L^{p'}_{\rm loc}$, $U$ is a matrix weight, and $W^{-1/p}h_j$ is locally integrable for each $h_j\in L^p(\mathbb C^d)$. In particular, for $Q\in\mathcal G$,
\begin{equation}\label{eq:domainlocal}
 \int_Q|W(y)^{-1/p}h_j(y)|\,dy\le\left(\int_Q\|W^{-1/p}\|_{\op}^{p'}\right)^{1/p'}\|h_j\|_{L^p(Q)}<\infty.
\end{equation}
The reducing and pointwise maximal operators are measurable and independent of the representatives almost everywhere. The norm of the countable sequence extension is the supremum of the finite-sequence norms.
\end{lemma}

\begin{lemma}[Dyadic reverse H\"older inequality]\label{lem:dyadicRH}
Let $1<s<\infty$, let $\D$ be a dyadic grid and let $v$ be a scalar weight in $A_s(\D)$ as defined in \eqref{eq:scalarAp}. Then
\begin{equation}\label{eq:dyadic-RH}
 \begin{aligned}
 1&\le [v]_{A_\infty(\D)}\le 2^{s+1}[v]_{A_s(\D)},\\
 \av{v^{1+\eta}}_Q&\le 2\av{v}_Q^{1+\eta}
 \quad\left(0<\eta\le\frac{1}{2^{n+2}[v]_{A_\infty(\D)}}\right).
 \end{aligned}
\end{equation}
The second inequality holds for every $Q\in\D$.
\end{lemma}

The dyadic formulation is needed for \Probref{prob:DKPS}. The corresponding estimate for all cubes is due to Hyt\"onen and P\'erez~\cite[Theorem~2.3]{HP}.

\begin{lemma}[Size and integrability of the matrix kernel]\label{lem:kernel}
Let $1<p<\infty$, let $\mathcal G=\mathcal Q$ or $\mathcal G=\D$ for a fixed dyadic grid, and suppose $W\in A_p(\mathcal G)$. Then there is a constant $0<c_n\le1$ depending only on $n$ such that, for $s=p'(1+c_n/[W^{-1/(p-1)}]_{A_{p',\infty}^{\mathrm{sc}}(\mathcal G)})$,
\begin{equation}\label{eq:kernel}
 \begin{aligned}
 1&\le [W^{-1/(p-1)}]_{A_{p',\infty}^{\mathrm{sc}}(\mathcal G)}
 \lesssim_{n,d,p}[W]_{A_p(\mathcal G)}^{1/(p-1)},\\
 \sup_{Q\in\mathcal G}\left(\fint_Q\|V_QW(y)^{-1/p}\|_{\op}^s\,dy\right)^{1/s}
 &\lesssim_{n,d,p}[W]_{A_p(\mathcal G)}^{1/p}.
 \end{aligned}
\end{equation}
\end{lemma}

\begin{lemma}[The norm of one reducing average]\label{lem:localnorm}
Let $W\in A_p(\mathcal G)$, $Q\in\mathcal G$ and $1\le q\le\infty$. Then the exact norm of this average is
\begin{equation}\label{eq:local-story}
 \|\overrightarrow B_Q\|_{L^p(\ell_d^q)\to L^p(\ell^q)}
 =\|B_Q\|_{L^p(\mathbb C^d)\to L^p}
 =\left(\fint_Q\|V_QW(y)^{-1/p}\|_{\op}^{p'}dy\right)^{1/p'}.
\end{equation}
The norm is independent of $q$.

For $W=wI_d$, take $V_Q=\av{w}_Q^{1/p}I_d$ and set $\sigma=w^{-1/(p-1)}$. Then the kernel is $(\av{w}_Q/w(y))^{1/p}$ and
\begin{equation}\label{eq:one-cube-cost}
 \|B_Q\|_{L^p(\mathbb C^d)\to L^p}=\left(\av{w}_Q\av{\sigma}_Q^{p-1}\right)^{1/p}.
\end{equation}
This is the local scalar $A_p$ characteristic to the power $1/p$. Equality is attained by $h=\sigma^{1/p}e\one_Q$ for any fixed unit vector $e$. Now let $W\in A_p(\mathcal G)$ be a general matrix weight. Then
\begin{equation}\label{eq:local-global-cost}
 \sup_{Q\in\mathcal G}\|\overrightarrow B_Q\|_{L^p(\ell_d^q)\to L^p(\ell^q)}\eqsim_{d,p} [W]_{A_p(\mathcal G)}^{1/p}.
\end{equation}
\end{lemma}

\begin{lemma}[Local sequence estimate]\label{lem:weak}
Let $Q$ be the root of a finite dyadic tree $\mathcal U$ and let $u_1,\ldots,u_N$ be nonnegative measurable functions on $Q$, with $N\in\mathbb N$. Suppose $1<q\le\infty$ and $F=\|(u_j)_j\|_{\ell^q}\in L^1(Q)$. Then for $\lambda>0$,
\begin{equation}\label{eq:weakfreeze}
 |\{x\in Q:\|(M_Q^{\mathcal U}u_j(x))_j\|_{\ell^q}>\lambda\}|
 \le (2^n+1)\frac{q'}{\lambda}\int_Q F.
\end{equation}
The constant is independent of the depth and of $N$.
\end{lemma}

\begin{proposition}[Sparse control of the sequence maximal function]\label{prop:freeze}
Let $W$ be a matrix weight, $1<p<\infty$, and fix a reducing family satisfying \eqref{eq:reducing} on a finite dyadic forest $\mathcal T$. Let $N\in\mathbb N$, $f_1,\ldots,f_N\in L^1_{\rm loc}(\R^n,\mathbb C^d)$ and $1<q\le\infty$. Then there is a $1/2$-sparse family
$\Spar\subseteq\mathcal T$ such that
\begin{equation}\label{eq:sparsefreeze}
 \left(\sum_{j=1}^{N}\left(\sup_{\substack{R\in\mathcal T\\x\in R}}
                \av{|V_Rf_j|}_R\right)^q\right)^{1/q}
 \lesssim_{n,d,p} q'
 \left(\sum_{Q\in\Spar}\one_Q(x)\left(\av{\|(V_Q f_j)_j\|_{\ell_{N,d}^q}}_Q\right)^q\right)^{1/q}.
\end{equation}
At $q=\infty$ the two sums raised to $1/q$ are interpreted as suprema of the corresponding entries. Each coordinate on the left has its own supremum. On the right the matrix is fixed and the sequence norm is inside the average. The bound holds for every $x\in\R^n$, with the empty maximum convention of \Subref{subsec:notation}. The implicit constant is independent of the weight, the depth of the trees and the number of coordinates.
\end{proposition}

\begin{lemma}[A scalar sparse estimate]\label{lem:scalar-sparse}
Let $1\le r<p<\infty$, $1\le q\le\infty$, and let $\Spar$ be a finite $1/2$-sparse subfamily of a dyadic grid $\D$. For $F\ge0$ in $L^p(\R^n)$, set $H=(M^\D(F^r))^{1/r}$. Then
\begin{equation}\label{eq:scalar-sparse}
 \left\|\left(\sum_{Q\in\Spar}\one_Q\av F_{r,Q}^{q}\right)^{1/q}\right\|_p
 \le C_p\|H\|_p,
\end{equation}
where the left side for $q=\infty$ is the norm of the supremum over $Q$. The constant is independent of $r,q$ and $\Spar$.
\end{lemma}

\begin{lemma}[Power weights]\label{lem:power}
Let $1<r<\infty$ and $-1<\alpha<r-1$. Then the weight $|x_1|^\alpha$ (assigned value one on $\{x_1=0\}$) satisfies
\begin{equation}\label{eq:powerformula}
 [|x_1|^\alpha]_{A_r}\eqsim_r
 \frac1{1+\alpha}\left(1-\frac{\alpha}{r-1}\right)^{1-r}.
\end{equation}
\end{lemma}

\subsection{The upper estimates}\label{subsec:upper}
We first explain why a direct estimate of each coordinate does not cover the full range of sequence exponents. On a fixed grid $\D$, suppose $W\in A_p(\D)$. To see the restriction, let $h_j\in L^p(\mathbb C^d)$, $\delta=c_n/[W^{-1/(p-1)}]_{A_{p',\infty}^{\mathrm{sc}}(\D)}$ and $r=p(1+\delta)/(1+p\delta)$, with $c_n$ from \Lref{lem:kernel}. The conjugate exponent of $r$ is $p'(1+\delta)$. H\"older's inequality and \eqref{eq:kernel} give
\begin{equation}\label{eq:coordinate-h}
 \widetilde{\mathcal M}_W^\D h_j\lesssim_{n,d,p}[W]_{A_p(\D)}^{1/p}\big(M^\D(|h_j|^r)\big)^{1/r}.
\end{equation}
To sum \eqref{eq:coordinate-h} by the ordinary vector maximal inequality, one would need both $p/r>1$ and $q/r>1$. For $1<q<p$,
\[
 q/r>1\quad\Longleftrightarrow\quad \delta>\frac{p-q}{p(q-1)}.
\]
As the dual directional characteristic grows, $r$ approaches $p$, and the right-hand condition can fail. We use a common sparse family before applying H\"older's inequality. Simply changing the order of the norm and the supremum does not achieve this: For nonnegative functions $u_j$ one has
\[
 \sup_{R\ni x}\|(\av{u_j}_R)_j\|_{\ell^q}
 \le\|(\sup_{R\ni x}\av{u_j}_R)_j\|_{\ell^q}.
\]
For finite $q$, the reverse inequality fails uniformly in $N$, even in dimension one. Take $u_j=\one_{[2^j,2^{j+1})}$, $1\le j\le N$, and $0<x<1$. Disjoint supports bound the left side by one. The standard dyadic intervals $[0,2^{j+1})$ give a right side at least $N^{1/q}/2$. The same argument works for all open intervals. At $q=\infty$ the two suprema commute.

For a fixed reducing matrix $V_Q$, the averages of $|V_QW^{-1/p}h_j|$ are ordinary scalar averages. The local weak estimate in \Lref{lem:weak} is the basis of the common sparse selection in \Pref{prop:freeze}. This is related to scalar sequence sparse domination in \cite{CLPR,CDO,HL}. The additional matrix stopping condition in \eqref{eq:badset} permits comparison of every descendant matrix with the selected ancestor matrix.

\begin{proof}[Proof of the upper estimates]
We first prove the grid estimate \eqref{eq:dyadic-main}, then \eqref{eq:main} and \eqref{eq:sharp-pure} for all cubes.

Fix a dyadic grid $\D$ and a matrix weight $W\in A_p(\D)$. For $1<q\le\infty$, take a finite sequence $\mathbf h=(h_j)_{j=1}^N$ in $L^p(\ell_{N,d}^q)$. Let $f_j=W^{-1/p}h_j$ for the use of \Pref{prop:freeze}. Set
\[
 F(x)=\left(\sum_{j=1}^N|h_j(x)|^q\right)^{1/q},\qquad \|F\|_p=\|\mathbf h\|_{L^p(\ell_{N,d}^q)}.
\]
At $q=\infty$, $F(x)=\max_{1\le j\le N}|h_j(x)|$. By \eqref{eq:domainlocal}, the functions $f_j$ are locally integrable. Define $\delta=c_n/[W^{-1/(p-1)}]_{A_{p',\infty}^{\mathrm{sc}}(\D)}$, $s=p'(1+\delta)$ and $r=s'$ as in \Lref{lem:kernel}. Then
\begin{equation}\label{eq:rexponent}
 r=\frac{p(1+\delta)}{1+p\delta},\qquad
 \frac{p}{r}=1+\frac{(p-1)\delta}{1+\delta},\qquad 1<r<p.
\end{equation}

Let $\mathcal T$ be a finite dyadic forest, and choose the sparse family $\Spar$ given by \Pref{prop:freeze}. For each $Q\in\Spar$ and almost every $y\in Q$, the matrix factor is common to the entire sequence:
\[
 \|(V_QW(y)^{-1/p}h_j(y))_j\|_{\ell_{N,d}^q}\le K_Q(y)\|(h_j(y))_j\|_{\ell_{N,d}^q}.
\]
Writing $K_Q(y)=\|V_QW(y)^{-1/p}\|_{\op}$, H\"older's inequality and \Lref{lem:kernel} give
\begin{equation}\label{eq:inside}
 \av{\|(V_Q f_j)_j\|_{\ell_{N,d}^q}}_Q
 \le\av{K_Q F}_Q\le\av{K_Q}_{s,Q}\av F_{r,Q}
 \lesssim_{n,d,p}[W]_{A_p(\D)}^{1/p}\av F_{r,Q}.
\end{equation}
Only $r<p$ is needed to estimate the scalar sparse sum.

By \eqref{eq:sparsefreeze} and \eqref{eq:inside}, for $1<q<\infty$ we obtain
\begin{equation}\label{eq:forest-global}
 \left\|\left(\sup_{\substack{R\in\mathcal T\\x\in R}}\av{|V_Rf_j|}_R\right)_j\right\|_{\ell^q}
 \lesssim_{n,d,p}q'[W]_{A_p(\D)}^{1/p}
 \left(\sum_{Q\in\Spar}\one_Q(x)\av F_{r,Q}^{q}\right)^{1/q}.
\end{equation}
At $q=\infty$, use a supremum on the right. By \Lref{lem:scalar-sparse}, the $L^p$ norm of the right side is at most $C_{n,d,p}q'[W]_{A_p(\D)}^{1/p}\|H\|_p$, where $H=(M^\D(F^r))^{1/r}$.
For the scalar maximal estimate we use $p/r>1$. By \eqref{eq:rexponent},
\[
 \frac{p}{p-r}=\frac{1+p\delta}{(p-1)\delta},\qquad
 \frac1r-\frac1p=\frac{(p-1)\delta}{p(1+\delta)}.
\]
Since $0<\delta\le1$, the bound $\delta\log(1/\delta)\le1/\exp(1)$ gives
\[
 \delta^{-(1/r-1/p)}\le\exp\big((p-1)/(p\exp(1))\big).
\]
Consequently,
\begin{align}\label{eq:gain}
 \|(M^\D(F^r))^{1/r}\|_p
 &\le \left(\frac{p}{p-r}\right)^{1/r}\|F^r\|_{p/r}^{1/r}\notag\\
 &\lesssim_p\delta^{-1/r}\|F\|_p
 \lesssim_p\delta^{-1/p}\|F\|_p
 \lesssim_{n,p}[W^{-1/(p-1)}]_{A_{p',\infty}^{\mathrm{sc}}(\D)}^{1/p}\|F\|_p.
\end{align}
Combining \eqref{eq:scalar-sparse} and \eqref{eq:gain}, we obtain
\begin{equation}\label{eq:sparseLp}
 \left\|\left(\sum_{Q\in\Spar}\one_Q\av F_{r,Q}^{q}\right)^{1/q}\right\|_p
 \lesssim_{n,p}[W^{-1/(p-1)}]_{A_{p',\infty}^{\mathrm{sc}}(\D)}^{1/p}\|F\|_p.
\end{equation}
This also holds for the sparse supremum when $q=\infty$. Together with \eqref{eq:forest-global}, it gives the desired bound uniformly on finite forests. Enumerate $\D$ as $R_1,R_2,\ldots$. The maximal members of $\{R_1,\ldots,R_k\}$ are disjoint. Completing their descendants down to the smallest side length in the collection gives a finite forest containing all these cubes. We obtain
\[
 \left\|\left(\max_{1\le i\le k}\big(\one_{R_i}(x)\av{|V_{R_i} f_j|}_{R_i}\big)\right)_j\right\|_{L^p(\ell^q)}
 \lesssim_{n,d,p}q'[W]_{A_p(\D)}^{1/p}[W^{-1/(p-1)}]_{A_{p',\infty}^{\mathrm{sc}}(\D)}^{1/p}\|F\|_p.
\]
The maxima increase to the dyadic suprema. Taking the $p$th powers and applying monotone convergence gives the dyadic sequence norm. Thus
\begin{equation}\label{eq:dyadic-allq}
 \|\bM_W^\D\|_{L^p(\ell_d^q)\to L^p(\ell^q)}
 \lesssim_{n,d,p}q'[W]_{A_p(\D)}^{1/p}[W^{-1/(p-1)}]_{A_{p',\infty}^{\mathrm{sc}}(\D)}^{1/p}.
\end{equation}
By \Lref{lem:domain}, the bound passes to countable sequences. Using \eqref{eq:kernel}, we have
\[
 [W]_{A_p(\D)}^{1/p}[W^{-1/(p-1)}]_{A_{p',\infty}^{\mathrm{sc}}(\D)}^{1/p}\lesssim_{n,d,p} [W]_{A_p(\D)}^{1/p+1/(p(p-1))}=[W]_{A_p(\D)}^{1/(p-1)}.
\]
This proves both estimates in \eqref{eq:dyadic-main}.

\smallskip
\noindent\emph{Passage from one grid to all cubes.}
We now assume the full-cube condition $W\in A_p$. For each of the $3^n$ vectors $u\in\{0,1/3,2/3\}^n$, let $\D^u$ be the dyadic grid with cubes
\[
 2^{-k}\big(m+(-1)^k u+[0,1)^n\big),\qquad k\in\mathbb Z,\quad m\in\mathbb Z^n.
\]
These are nested grids because $3u\in\mathbb Z^n$. Given an open cube $Q$, choose a dyadic length $L$ such that $3\ell(Q)<L\le6\ell(Q)$. At scale $L$ the three boundary sets in each coordinate are separated by $L/3$. Each coordinate projection of $Q$ has length less than $L/3$, so one may choose a shift avoiding that projection. This gives a cube $R\in\D^u$ with $Q\subseteq R$ and $|R|/|Q|\le 6^n$. By \eqref{eq:reducing},
\begin{align}\label{eq:adj}
 |V_Qv|&\le c_d\left(\frac{|R|}{|Q|}\right)^{1/p}|V_Rv|,\notag\\
 \av{|V_Q f_j|}_Q&\le c_d\left(\frac{|R|}{|Q|}\right)^{1+1/p}\av{|V_R f_j|}_R.
\end{align}
Taking the supremum over $Q$ in \eqref{eq:adj}, we obtain
\[
 \widetilde{\mathcal M}_W h_j\le c_d6^{n(1+1/p)}\sum_u\widetilde{\mathcal M}_W^{\D^u}h_j.
\]
By the triangle inequalities in $\ell^q$ and $L^p$ and \eqref{eq:dyadic-allq},
\begin{align*}
 \|\bM_W\|_{L^p(\ell_d^q)\to L^p(\ell^q)}
 &\lesssim_{n,d,p}\sum_u\|\bM_W^{\D^u}\|_{L^p(\ell_d^q)\to L^p(\ell^q)}\\
 &\lesssim_{n,d,p}q'[W]_{A_p}^{1/p}[W^{-1/(p-1)}]_{A_{p',\infty}^{\mathrm{sc}}}^{1/p}\\
 &\lesssim_{n,d,p}q'[W]_{A_p}^{1/(p-1)}.
\end{align*}
The second line uses restriction of the characteristics to $\D^u$ and the last uses \eqref{eq:kernel}.

\end{proof}

\subsection{Sharpness and the mixed characterization}\label{sub:tests}
By \eqref{eq:local-global-cost} and $B_Qh\le\widetilde{\mathcal M}_Wh$,
\[
 [W]_{A_p}^{1/p}\lesssim_{d,p}\|\bM_W\|_{L^p(\ell_d^q)\to L^p(\ell^q)}.
\]
The first power-weight family has uniformly bounded dual directional characteristic. The second saturates the upper comparison in \eqref{eq:kernel}. We use one sequence entry throughout.

For $W=wI_d$, the characteristic \eqref{eq:A} equals \eqref{eq:scalarAp} and the dual directional characteristic is $[w^{-1/(p-1)}]_{A_\infty}$. Write $x=(x_1,x')\in\R\times\R^{n-1}$ and use weights depending only on $x_1$. When $n=1$, the transverse variable is absent and its volume factor is one. Power weights are assigned the value one on the null set $\{x_1=0\}$. Their cube averages equal interval averages in the first coordinate.

\begin{proof}[Proof of sharpness and sufficiency]
\noindent\emph{Necessity of $a\ge1/p$.} Take $0<\varepsilon<1/2$ and $w_-(x)=|x_1|^{-1+\varepsilon}$ off $\{x_1=0\}$, define $W_-=w_-I_d$. By \eqref{eq:powerformula},
\[
 [W_-]_{A_p}\eqsim_p\frac{1}{\varepsilon}\left(1+\frac{1-\varepsilon}{p-1}\right)^{1-p}\eqsim_p\varepsilon^{-1},\qquad
 w_-^{-1/(p-1)}=|x_1|^{(1-\varepsilon)/(p-1)}.
\]
On every interval $I$, the set $\{u\in I:|u|\ge\frac12\sup_{y\in I}|y|\}$ has measure at least $|I|/4$, by considering the endpoint of largest absolute value. Thus for $\sigma=w_-^{-1/(p-1)}$ and any cube $Q$,
\[
 \operatorname*{ess\,sup}_{y\in Q}\sigma(y)\le 4\cdot 2^{1/(p-1)}\av{\sigma}_Q,\qquad
 [w_-^{-1/(p-1)}]_{A_\infty}=\sup_Q\frac{\int_Q M(\sigma\one_Q)}{\sigma(Q)}
 \le 4\cdot 2^{1/(p-1)} .
\]
Choose a unit vector $e\in\mathbb C^d$ and take the single-entry sequence given by $h=w_-^{1/p}e\one_{(1/2,1)\times(0,1)^{n-1}}$. On this support $1\le w_-\le2$, so $1/2\le\|h\|_p^p\le1$. Testing on $Q=(0,1)^n$ and using the lower inequality in \eqref{eq:reducing} gives
\begin{equation}\label{eq:minus}
 \widetilde{\mathcal M}_{W_-}h(x)\ge\tfrac12\av{w_-}_Q^{1/p}
 =\tfrac12\varepsilon^{-1/p}\quad (x\in Q).
\end{equation}
Since $\|h\|_p\le1$, \eqref{eq:minus} and \eqref{eq:admissible} imply
\[
 \varepsilon^{-1/p}\lesssim_p\|\bM_{W_-}\|_{L^p(\ell_d^q)\to L^p(\ell^q)}\lesssim\varepsilon^{-a}.
\]
Letting $\varepsilon\downarrow0$ gives $a\ge1/p$.

\noindent\emph{Necessity of $(p-1)a+b\ge1$.} Let $0<\varepsilon<1/2$, $\alpha=(p-1)(1-\varepsilon)$ and $w_+(x)=|x_1|^\alpha$. Define $W_+=w_+I_d$. Then
\[
 [W_+]_{A_p}\eqsim_p\varepsilon^{-(p-1)},\qquad [w_+^{-1/(p-1)}]_{A_\infty}\lesssim_{n,p}\varepsilon^{-1}.
\]
Here $1-\alpha/(p-1)=\varepsilon$ and $1+\alpha\eqsim_p1$. The upper bound for the dual characteristic follows from \eqref{eq:kernel}. Set $Q_0=[0,1)^n$, $\sigma=w_+^{-1/(p-1)}$ and $h=\sigma^{1/p}e\one_{Q_0}$ for a unit vector $e$. Then $W_+^{-1/p}h=\sigma e\one_{Q_0}$ and $\|h\|_p^p=\sigma(Q_0)=\varepsilon^{-1}$. For almost every $x\in Q_0$, choose a dyadic length $L$ with $x_1<L\le2x_1$, $L\le1$. Let $R\subseteq Q_0$ be the dyadic cube of side $L$ containing $x$ whose first interval is $[0,L)$. Its interior is also an admissible open cube for almost every such $x$. Direct integration gives
\[
 \av{\sigma}_R=\varepsilon^{-1}L^{-1+\varepsilon}
 \ge(2\varepsilon)^{-1}\sigma(x),\qquad
 \av{w_+}_R=L^\alpha/(1+\alpha).
\]
Since $M(\sigma\one_{Q_0})(x)\ge\av{\sigma}_R$ almost everywhere on $Q_0$, we obtain
\[
 [w_+^{-1/(p-1)}]_{A_\infty}\ge\frac{1}{\sigma(Q_0)}\int_{Q_0}M(\sigma\one_{Q_0})\ge(2\varepsilon)^{-1}.
\]
The lower reducing inequality on $R$ gives, for both the full operator and its standard-grid restriction,
\[
 \av{|V_RW_+^{-1/p}h|}_R
 \ge \av{w_+}_R^{1/p}\av{\sigma}_R
 = \frac{L^{(-1+\varepsilon)/p}}{\varepsilon(1+\alpha)^{1/p}}
 \gtrsim_p \varepsilon^{-1}|h(x)|.
\]
Using the interior of $R$ gives the same lower estimate for any admissible full-cube reducing family. Thus
\begin{equation}\label{eq:plus}
 [w_+^{-1/(p-1)}]_{A_\infty}\eqsim_{n,p}\varepsilon^{-1},\qquad
 \|\bM_{W_+}\|_{L^p(\ell_d^q)\to L^p(\ell^q)}
 \gtrsim_p\varepsilon^{-1}.
\end{equation}
If \eqref{eq:admissible} holds, then
\[
 \varepsilon^{-1}\lesssim_p\|\bM_{W_+}\|_{L^p(\ell_d^q)\to L^p(\ell^q)}\lesssim\varepsilon^{-((p-1)a+b)}.
\]
Letting $\varepsilon\downarrow0$ gives $(p-1)a+b\ge1$.

\noindent\emph{The optimal pure norm law.} For the lower bound in \eqref{eq:extremal-norm}, the preceding estimates give a constant $C_0\ge1$, depending only on $p$, such that $[W_+]_{A_p}\le C_0\varepsilon^{-(p-1)}$ for $0<\varepsilon<1/2$. Given $t>2^{p-1}C_0$, choose $\varepsilon=(C_0/t)^{1/(p-1)}$. Then $[W_+]_{A_p}\le t$, and \eqref{eq:plus} gives
\[
 \sup_{\substack{W\in A_p\\ {[W]}_{A_p}\le t}}\|\bM_W\|_{L^p(\ell_d^q)\to L^p(\ell^q)}
 \gtrsim_{n,p} t^{1/(p-1)}.
\]
For $1\le t\le 2^{p-1}C_0$, the constant weight $I_d$ belongs to the allowed class and gives norm at least one, by the test of a single entry equal to a unit vector on one cube. This yields the same lower bound after changing the constant. The upper bound in \eqref{eq:extremal-norm} follows from \eqref{eq:sharp-pure}. All the lower estimates use the lower reducing inequality in \eqref{eq:reducing}, so they hold regardless of the chosen reducing families.

\noindent\emph{Sharpness on every fixed grid.}
Choose any $Q_*\in\D$, with lower corner $z$ and side length $L_*>0$. Use the change of variables $x=z+L_*u$. Define
\[
 W_{\pm,*}(x)=W_\pm\big((x-z)/L_*\big).
\]
The full-cube scalar characteristics are unchanged by this translation and dilation. The descendants of $Q_*$ correspond exactly to the standard dyadic descendants of $[0,1)^n$. Testing on $Q_*$ and using the full-cube upper bounds gives
\begin{align*}
 [W_{-,*}]_{A_p(\D)}&\eqsim_p\varepsilon^{-1},&
 [W_{-,*}^{-1/(p-1)}]_{A_{p',\infty}^{\mathrm{sc}}(\D)}&\eqsim_{n,p}1,\\
 [W_{+,*}]_{A_p(\D)}&\eqsim_p\varepsilon^{-(p-1)},&
 [W_{+,*}^{-1/(p-1)}]_{A_{p',\infty}^{\mathrm{sc}}(\D)}&\eqsim_{n,p}\varepsilon^{-1}.
\end{align*}
For the last lower bound, the descendants of $Q_*$ give the same estimate for the dual weight as the descendants of $Q_0=[0,1)^n$ in the positive-power test. For either sign, let $h$ be the corresponding single-entry test above and take
\[
 h_*(x)=L_*^{-n/p}h\big((x-z)/L_*\big),\qquad \|h_*\|_p=\|h\|_p.
\]
For a dyadic descendant $R$ of $[0,1)^n$, write $R_*=z+L_*R$. The reducing inequality on $R_*$ and the fact that the weight is scalar times $I_d$ give
\[
 \av{|V_{R_*}W_{\pm,*}^{-1/p}h_*|}_{R_*}
 \ge L_*^{-n/p}\av{w_\pm}_R^{1/p}\av{|w_\pm^{-1/p}h|}_R.
\]
In the normalized cube, these are the averages used in \eqref{eq:minus} and \eqref{eq:plus}. The change of variables cancels $L_*^{-n/p}$ in the $L^p$ norm. Thus
\begin{align*}
 \|\bM_{W_{-,*}}^\D\|_{L^p(\ell_d^q)\to L^p(\ell^q)}&\gtrsim_p\varepsilon^{-1/p},\\
 \|\bM_{W_{+,*}}^\D\|_{L^p(\ell_d^q)\to L^p(\ell^q)}&\gtrsim_p\varepsilon^{-1}.
\end{align*}
Both necessary conditions in \eqref{eq:region} follow. To obtain the lower estimate in \eqref{eq:dyadic-extremal}, take $\varepsilon=(C_0/t)^{1/(p-1)}$ when $t>2^{p-1}C_0$, where $C_0$ is as in the full-cube test. Then
\[
 [W_{+,*}]_{A_p(\D)}\le t,\qquad
 \|\bM_{W_{+,*}}^\D\|_{L^p(\ell_d^q)\to L^p(\ell^q)}\gtrsim_{n,p}t^{1/(p-1)}.
\]
For $1\le t\le2^{p-1}C_0$, use $W=I_d$ and one coordinate. The upper estimate follows from \eqref{eq:dyadic-main}. All constants are independent of $Q_*$ and $\D$. Sufficiency for the mixed region follows below.

\noindent\emph{Sufficiency.} For any fixed unit vector $e$,
\[
 \frac{|W(x)^{1/p}e|}{|W(y)^{1/p}e|}\le\|W(x)^{1/p}W(y)^{-1/p}\|_{\op}.
\]
Thus $1\le[|W^{1/p}e|^p]_{A_p}\le [W]_{A_p}$. Also $[W^{-1/(p-1)}]_{A_{p',\infty}^{\mathrm{sc}}}\ge1$ by differentiation. If $b\ge1/p$ and $a\ge1/p$, the mixed estimate gives \eqref{eq:admissible}. When $0\le b<1/p$, \eqref{eq:kernel} yields
\begin{align*}
 [W]_{A_p}^{1/p}[W^{-1/(p-1)}]_{A_{p',\infty}^{\mathrm{sc}}}^{1/p}
 &\lesssim_{n,d,p}[W]_{A_p}^{1/p+(1/p-b)/(p-1)}[W^{-1/(p-1)}]_{A_{p',\infty}^{\mathrm{sc}}}^b\\
 &=[W]_{A_p}^{(1-b)/(p-1)}[W^{-1/(p-1)}]_{A_{p',\infty}^{\mathrm{sc}}}^b.
\end{align*}
The condition $(p-1)a+b\ge1$ gives $a\ge(1-b)/(p-1)$. Since $[W]_{A_p}\ge1$, this proves sufficiency in the remaining range. The same reasoning with the restricted characteristics proves the asserted dyadic sufficiency. This concludes the mixed characterization.
\end{proof}

\begin{proof}[Proof of \Cref{cor:tradeoff}]
Solving the two inequalities in \eqref{eq:region} for $a$ gives \eqref{eq:tradeoff}. For $0\le\theta\le1$,
\begin{align*}
 [W]_{A_p}^{1/p}[W^{-1/(p-1)}]_{A_{p',\infty}^{\mathrm{sc}}}^{1/p}
 &\lesssim_{n,d,p}[W]_{A_p}^{1/p+(1-\theta)/(p(p-1))}[W^{-1/(p-1)}]_{A_{p',\infty}^{\mathrm{sc}}}^{\theta/p}.
\end{align*}
Use \eqref{eq:main} and the identity $1/p+(1-\theta)/(p(p-1))=(1-\theta)/(p-1)+\theta/p$ to obtain \eqref{eq:frontier}. Every pair on this segment satisfies $(p-1)a+b=1$, so a strict reduction in either admissible coordinate violates \eqref{eq:region}.
\end{proof}

\subsection{The finite-sequence estimate}\label{subsec:length-proof}
\begin{proof}[Proof of \Cref{cor:length}]
Set $L=1+\log N$ and $s=1+1/L$. For $1\le q<s$, the finite-dimensional norm comparison is
\[
 \|u\|_{\ell_N^s}\le\|u\|_{\ell_N^q}\le N^{1/q-1/s}\|u\|_{\ell_N^s},
 \qquad N^{1/q-1/s}\le\exp(1).
\]
Consequently,
\[
 \|\bM_{W,N}^{\mathcal G}\|_{L^p(\ell_{N,d}^q)\to L^p(\ell_N^q)}
 \le\exp(1)\|\bM_{W,N}^{\mathcal G}\|_{L^p(\ell_{N,d}^s)\to L^p(\ell_N^s)}.
\]
For $\mathcal G=\mathcal Q$ use \eqref{eq:main}, and for $\mathcal G=\D$ use \eqref{eq:dyadic-main}, at exponent $s$ if $q<s$ and at $q$ otherwise. Since $s'=1+L\le2L$, the bound \eqref{eq:length} follows from
\[
 q'>s'\ge L\quad(q<s),\qquad q'\le s'\le2L\quad(q\ge s).
\]

For $W=I_d$ and $V_Q=I_d$, the upper bound in \eqref{eq:length-sharp} is the scalar case of \eqref{eq:length}. To prove the lower bound, take any cube $Q$ of the given grid (the standard grid if $\mathcal G=\mathcal Q$). Let $m=\lfloor\log_2 N\rfloor$. For $N\ge2$, partition $Q$ into $2^m$ slabs $E_j$ of equal width in its first coordinate, and set $h_j=e\one_{E_j}$, $1\le j\le 2^m$, for a fixed unit vector $e$. Pad the sequence with zeros to length $N$. Then
\[
 \|(h_j)_j\|_{L^p(\ell_{N,d}^q)}=|Q|^{1/p}.
\]
For almost every $x\in Q$, let $Q_\nu(x)$, $0\le\nu\le m$, be the unique descendant of $Q$ containing $x$ and having side length $2^{\nu-m}\ell(Q)$. The first-coordinate projection of $Q_\nu(x)$ contains exactly $2^\nu$ of the slab intervals. For $1\le\nu\le m$, exactly $2^{\nu-1}$ of these intervals lie outside the projection of $Q_{\nu-1}(x)$. For their distinct indices $j$,
\[
 \av{|h_j|}_{Q_\nu(x)}=2^{-\nu},\qquad
 \widetilde{\mathcal M}_{I_d}^{\mathcal G}h_j(x)\ge2^{-\nu}.
\]
For $\mathcal G=\mathcal Q$, the interior of the same cube gives the average for almost every $x$. Thus, for $1\le q<\infty$,
\begin{equation}\label{eq:harmonic}
 \|\bM_{I_d,N}^{\mathcal G}\|_{L^p(\ell_{N,d}^q)\to L^p(\ell_N^q)}
 \ge \left(\frac12\sum_{\nu=1}^m 2^{-\nu(q-1)}\right)^{1/q}.
\end{equation}
At $q=1$ this is $m/2$. If $1<q\le2$, write $\tau=q-1$. Then
\[
 \sum_{\nu=1}^m 2^{-\nu\tau}=\frac{1-2^{-m\tau}}{2^\tau-1}\eqsim \min\{m,1/\tau\}.
\]
For $v=\min\{m,1/\tau\}\ge1$, we have
\[
 v^{1/q}\ge v\exp(-\tau\log(1/\tau))\ge v\exp(-1/\exp(1)).
\]
Since $m\eqsim L$ for $N\ge2$ and $q'\eqsim1/\tau$ in this range, \eqref{eq:harmonic} proves the required lower bound. For $2\le q\le\infty$ or $N=1$, one nonzero entry gives norm at least one. Finally, finite sequences extend by zero, so letting $N\to\infty$ proves the countable assertions.
\end{proof}

\subsection{The mixed square-maximal estimate}\label{sub:mixed-application}
Let $\D^1,\D^2$ be dyadic grids on $\R$ and $\mathcal R=\{I\times J:I\in\D^1,J\in\D^2\}$. For a matrix weight $W$ on $\R^2$ let $[W]_{A_p(\mathcal R)}$ be the characteristic defined by the right side of \eqref{eq:A} with the supremum over rectangles $R\in\mathcal R$ and both integrals over $R$. Assume $[W]_{A_p(\mathcal R)}<\infty$. For $I\in\D^1$, let $\psi_I=|I|^{-1/2}(\one_{I_-}-\one_{I_+})$, where $I_-$ and $I_+$ are its left and right halves. For $h\in L^p(\R^2,\mathbb C^d)$ define the Haar coefficients of $W^{-1/p}h$ by
\[
 f_I(x_2)=\int_\R W(y_1,x_2)^{-1/p}h(y_1,x_2)\psi_I(y_1)\,dy_1.
\]
Write $W_{x_1}(y_2)=W(x_1,y_2)$ and let $V_{x_1,J}$ be a reducing matrix for this sliced weight at exponent $p$ on $J$. We choose the sliced reducing matrices measurably as in \cite[Section~3.2 and Proposition~9.8]{DKPS}. Define
\begin{equation}\label{eq:mixed-operator-app}
 \mathcal{SM}_{W}h(x_1,x_2)=\left(\sum_{I\in\D^1}\frac{\one_I(x_1)}{|I|}
 \left(\sup_{\substack{J\in\D^2\\x_2\in J}}|V_{x_1,J}\av{f_I}_J|\right)^2\right)^{1/2}.
\end{equation}
The weight is incorporated before taking the Haar coefficients, and the vector average is taken before the Euclidean norm. Define also the corresponding first-variable square function by
\[
 \mathcal S_W^{(1)}h(x_1,x_2)=\left(\sum_{I\in\D^1}\frac{\one_I(x_1)}{|I|}|W(x_1,x_2)^{1/p}f_I(x_2)|^2\right)^{1/2}.
\]
The exponent in \cite[(5.2)--(5.3)]{DKPS} is
\begin{equation}\label{eq:gamma}
 \gamma(p)=\begin{cases}1/(p-1),&1<p\le2,\\[2pt]1/2+1/(p(p-1)),&p>2.\end{cases}
\end{equation}

\begin{remark}[A mixed square-maximal estimate]\label{rem:application}
Let $1<p<\infty$ and let $W$ be a matrix weight on $\R^2$ satisfying $[W]_{A_p(\mathcal R)}<\infty$ for the product grid $\mathcal R$ defined above. Then, for the mixed operator $\mathcal{SM}_W$ in \eqref{eq:mixed-operator-app} and the exponent $\gamma(p)$ in \eqref{eq:gamma},
\begin{equation}\label{eq:mixed-app-bound}
 \|\mathcal{SM}_W\|_{L^p(\R^2,\mathbb C^d)\to L^p(\R^2)}
 \lesssim_{d,p} [W]_{A_p(\mathcal R)}^{\gamma(p)+1/(p-1)},\qquad 1<p<\infty.
\end{equation}
The same estimate holds with the two variables exchanged. This extends the range of \cite[Lemma~8.2(1)]{DKPS} from $1<p\le2$ to $1<p<\infty$ for the indicated mixed operator.
\end{remark}

\begin{proof}[Proof of \Remref{rem:application}]
The argument of \Lref{lem:domain} applied to rectangles gives local integrability of $W^{-1/p}$ in exponent $p'$, hence of $f=W^{-1/p}h$ for $h\in L^p(\R^2,\mathbb C^d)$. The grids are countable. By Fubini, the coefficient functions are defined almost everywhere after modification on a common null set in the second variable. Moreover,
\[
 \int_J|f_I(x_2)|\,dx_2\le |I|^{-1/2}\int_{I\times J}|f|<\infty.
\]
Each $f_I$ is locally integrable, since a bounded interval is covered by finitely many grid intervals of a fixed length. We may first take the sum over a finite family $\mathcal I\subseteq\D^1$. Fubini also gives local integrability of the sliced weights and their dual weights, for almost every fixed coordinate. By \cite[Lemma~3.6 with the dual weight]{DKPS}, on a set of full measure in each coordinate we have
\[
 [W(x_1,\cdot)]_{A_p(\D^2)},\quad [W(\cdot,x_2)]_{A_p(\D^1)}\lesssim_{d,p} [W]_{A_p(\mathcal R)}.
\]
For each such $x_1$ and $I\in\mathcal I$, define
\[
 g_I(x_2)=|I|^{-1/2}\one_I(x_1)W(x_1,x_2)^{1/p}f_I(x_2).
\]
Apply the known one-variable square-function bound \cite[(5.2)--(5.3)]{DKPS} on the slices $W(\cdot,x_2)$ for almost every $x_2$:
\[
 \|\mathcal S_W^{(1)}h(\cdot,x_2)\|_{L^p(dx_1)}
 \lesssim_{d,p} [W]_{A_p(\mathcal R)}^{\gamma(p)}
 \|h(\cdot,x_2)\|_{L^p(dx_1)}.
\]
Integrating its $p$th power in $x_2$ gives
\begin{equation}\label{eq:square-app}
 \|\mathcal S_W^{(1)}h\|_{L^p(\R^2)}\lesssim_{d,p}[W]_{A_p(\mathcal R)}^{\gamma(p)}\|h\|_{L^p(\R^2,\mathbb C^d)}.
\end{equation}
For almost every $x_1$, $(g_I)_I$ belongs to the unweighted sequence domain in \Cref{cor:dyadic}. The identity $W_{x_1}^{-1/p}g_I=|I|^{-1/2}\one_I(x_1)f_I$ and the triangle inequality under the integral give
\begin{align*}
 &\left\|\left(\sum_{I\in\mathcal I}\frac{\one_I(x_1)}{|I|}
 \sup_{\substack{J\in\D^2\\x_2\in J}}|V_{x_1,J}\av{f_I}_J|^2\right)^{1/2}\right\|_{L^p(dx_2)}\\
 &\quad\le\left\|\left(\sum_{I\in\mathcal I}(\widetilde{\mathcal M}_{W_{x_1}}^{\D^2}g_I)^2\right)^{1/2}\right\|_{L^p(dx_2)}\\
 &\quad\lesssim_{d,p} [W]_{A_p(\mathcal R)}^{1/(p-1)}\|\mathcal S_W^{(1)}h(x_1,\cdot)\|_{L^p(dx_2)}.
\end{align*}
Integrate the $p$th powers in $x_1$. By \eqref{eq:square-app}, the resulting bound is
\[
 C_{d,p}[W]_{A_p(\mathcal R)}^{p(\gamma(p)+1/(p-1))}\|h\|_{L^p(\R^2,\mathbb C^d)}^p.
\]
Let the finite families exhaust $\D^1$. All sums increase, and measurability follows from the sliced choices and the countable grids. Monotone convergence proves \eqref{eq:mixed-app-bound}. Exchanging $x_1$ and $x_2$ proves the symmetric estimate.
\end{proof}

\subsection{The estimate inherited from the single-function result}\label{subsec:comparison-proof}
\begin{proof}[Justification of the single-function comparison]
Let $W\in A_p$ and write $f=W^{-1/p}h$ for $h\in L^p(\mathbb C^d)$. By \eqref{eq:product}, the dual reducing matrix $V'_Q$ satisfies
\[
 \av{|V_Q f|}_Q\le\|V_QV'_Q\|_{\op}\av{|(V'_Q)^{-1}f|}_Q
 \lesssim_{d,p}[W]_{A_p}^{1/p}\av{|(V'_Q)^{-1}f|}_Q.
\]
We apply the published single-function mixed estimate of Kakaroumpas and Nieraeth~\cite{KN} with one factor, averaging exponent one and multiplier weight $W^{1/p}$. In their convex-set-valued formulation, a vector $f(y)$ can be represented by the balanced disk $\{zf(y):z\in\mathbb C,\ |z|\le1\}$, whose radius is $|f(y)|$. Applying a matrix to this disk gives the radius of the transformed vector. Their reducing matrix for the inverse multiplier at exponent $p'$ is comparable to $V'_Q$ in the sense of \eqref{eq:reducing}. Their local directional Fujii--Wilson characteristic for that multiplier is at most $[W^{-1/(p-1)}]_{A_{p',\infty}^{\mathrm{sc}}}^{1/p'}$, because the local dyadic maximal operator on a cube is majorized by the full scalar maximal operator. The exponent $p'/p$ in their estimate therefore gives
\[
 \left\|\sup_{Q\ni\,\cdot}\av{|(V'_Q)^{-1}f|}_Q\right\|_p
 \lesssim_{n,d,p} [W^{-1/(p-1)}]_{A_{p',\infty}^{\mathrm{sc}}}^{1/p}\|h\|_{L^p(\mathbb C^d)}.
\]
Together with \eqref{eq:product}, this gives the single-function case of \eqref{eq:main}. Write $K_Q(y)=\|V_QW(y)^{-1/p}\|_{\op}$ and, for a nonnegative scalar function $H$, set $TH=\sup_{Q\ni\,\cdot}\fint_QK_QH$. The column inequality gives
\[
 TH\le\sum_{k=1}^d\widetilde{\mathcal M}_W(He_k),\qquad
 \|T\|_{L^p\to L^p}\le d\|\widetilde{\mathcal M}_W\|_{L^p(\mathbb C^d)\to L^p}.
\]
For a finite sequence $\mathbf h$, let $v_j=|h_j|$ and $F=\|(v_j)_j\|_{\ell^q}$. At $q=p$ the bound follows by summation, as in \eqref{eq:diagonal-norm}, and at $q=\infty$ from $\sup_j Tv_j\le T(\sup_jv_j)$. For $p<q<\infty$ set $\theta=p/q$ and $u_j=F(v_j/F)^{1/\theta}$ on $\{F>0\}$, with value zero otherwise. Then $\sum_j u_j^p=F^p$. H\"older's inequality for the nonnegative measure $K_Q(y)\,dy/|Q|$ gives $Tv_j\le (Tu_j)^\theta(TF)^{1-\theta}$. Consequently
\[
 \|(Tv_j)_j\|_{L^p(\ell^q)}
 \le\left(\sum_j\|Tu_j\|_p^p\right)^{\theta/p}\|TF\|_p^{1-\theta}
 \le\|T\|_{L^p\to L^p}\|F\|_p.
\]
Since $\widetilde{\mathcal M}_W h_j\le Tv_j$, the mixed bound follows for $q\ge p$. Countable sequences follow by \Lref{lem:domain}.
\end{proof}

\Needspace{8\baselineskip}
\section{Proofs of the auxiliary lemmas}\label{sec:aux}
\subsection{The domain and the reverse H\"older inequality}\label{subsec:domain-proof}
\begin{proof}[Proof of \Lref{lem:domain}]
Set $W=I_d$ on an exceptional null set if necessary. If $[W]_{A_p(\mathcal G)}<\infty$, choose for each $Q\in\mathcal G$ an $x_0\in Q$ such that the inner integral in \eqref{eq:A} is finite. From
\[
 \|W(y)^{-1/p}\|_{\op}\le\|W(x_0)^{-1/p}\|_{\op}\|W(x_0)^{1/p}W(y)^{-1/p}\|_{\op}
\]
we obtain $\|W^{-1/p}\|_{\op}\in L^{p'}(Q)$. Consequently, for every $h_j\in L^p(\mathbb C^d)$,
\[
 \int_Q|W(y)^{-1/p}h_j(y)|\,dy\le\left(\int_Q\|W^{-1/p}\|_{\op}^{p'}\right)^{1/p'}\|h_j\|_{L^p(Q)}<\infty.
\]
When $\mathcal G$ is a dyadic grid, each bounded cube is covered by finitely many grid cubes of one fixed generation. Hence $W^{-1/p}h_j$ is locally integrable in both settings. The dual matrix weight is $U=W^{-1/(p-1)}$. Its local integrability follows from $\|U\|_{\op}=\|W^{-1/p}\|_{\op}^{p'}$. In particular, $U$ satisfies the definition of a matrix weight even before its $A_{p'}$ characteristic is estimated.

All individual averages in \eqref{eq:max} and \eqref{eq:pointwise-max} are finite by \eqref{eq:domainlocal}.
The full reducing operator is measurable: its strict superlevel sets are unions of open cubes for which the defining average exceeds the threshold. For the pointwise operator, we may restrict to cubes with rational lower corners and rational side lengths. Each open cube containing $x$ admits approximating cubes of this form containing it, with convergent endpoints. The corresponding averages converge since $W^{-1/p}h$ is locally integrable and $W(x)$ is fixed. This gives a countable supremum of measurable functions. Dyadic suprema are measurable by countability. Changing the representative of $h$ does not affect either maximal operator, and changing $W$ on a null set does not affect the operators almost everywhere, with the same reducing family.

For countable sequences, truncation to the first $N$ entries increases the $p$th powers of both sequence lengths to their countable versions. Monotone convergence shows that the operator norm of the countable extension equals the supremum of the finite-sequence norms. This includes $q=\infty$ and the restriction to a single grid.
\end{proof}

\begin{proof}[Proof of \Lref{lem:dyadicRH}]
Write $F_Q(x)=M_Q^\D v(x)$ on $Q$. For an integer $k\ge0$, let $F_Q^{(k)}$ be the local maximal function restricted to $Q$ and its descendants with side length at least $2^{-k}\ell(Q)$.
For $x\in Q$, $M^\D(v\one_Q)(x)=M_Q^\D v(x)$, because the average over an ancestor of $Q$ is at most $\av{v}_Q$. Dyadic differentiation gives $F_Q\ge v$ almost everywhere.

First, H\"older's inequality on a measurable $E\subseteq P$ with $P\in\D$ gives
\begin{align*}
 |E|^s &\le v(E)\left(\int_P v^{-1/(s-1)}\right)^{s-1},\\
 \frac{v(E)}{v(P)}&\ge\frac{1}{[v]_{A_s(\D)}}\left(\frac{|E|}{|P|}\right)^s.
\end{align*}
Starting at $Q$, stop on maximal dyadic subcubes $P'$ of each current cube $P$ for which $\av{v}_{P'}>2\av{v}_P$. Work first with finitely many generations. The stopped children occupy at most half of $P$. Let $E_P$ be $P$ minus the union of the cubes selected immediately from $P$. If $x\in E_P$, every descendant of $P$ containing $x$ has average at most $2\av{v}_P$. The same bound holds for ancestors of $P$ in $Q$: an intermediate ancestor between $P$ and its current parent has average at most twice the average of that parent by maximality of $P$, while each earlier selected ancestor has smaller average than its selected child by the stopping rule. Thus the truncated maximal function is at most $2\av{v}_P$ on $E_P$. The sets $E_P$ partition $Q$ in the finite tree. Hence
\begin{align*}
 \int_Q F_Q^{(k)}
 &\le 2\sum_P\av{v}_P |E_P|\le 2\sum_P v(P)\\
 &\le 2^{s+1}[v]_{A_s(\D)}\sum_P v(E_P)=2^{s+1}[v]_{A_s(\D)}v(Q),
\end{align*}
where $F_Q^{(k)}$ is the maximal function over the first $k$ generations. Monotone convergence proves the first bound in \eqref{eq:dyadic-RH}.

Write $a=\av{v}_Q$. For $\lambda\ge a$, let $\{P_i\}$ be the maximal strict subcubes of $Q$ with $\av{v}_{P_i}>\lambda$. Up to null sets these are disjoint and cover $\{F_Q>\lambda\}$. Each parent has average at most $\lambda$, so $v(P_i)\le2^n\lambda|P_i|$. On $P_i$ the maximal function $F_Q$ equals $F_{P_i}$, since every proper ancestor in $Q$ has average at most $\lambda$. Consequently
\begin{align}
 \int_{\{F_Q>\lambda\}}F_Q
 &\le [v]_{A_\infty(\D)}\sum_i v(P_i)
 \le 2^n [v]_{A_\infty(\D)}\lambda|\{F_Q>\lambda\}|,\label{eq:RHdist1}\\
 \int_{\{F_Q>\lambda\}}v
 &\le 2^n\lambda|\{F_Q>\lambda\}|.\label{eq:RHdist2}
\end{align}
For $H_T=\min(F_Q,T)$, $T>a$, the layer-cake formula applied to $H_T^\eta$ yields
\begin{align*}
 \int_Q H_T^{1+\eta}
 &=a^\eta\int_Q H_T+\eta\int_a^T\lambda^{\eta-1}
                         \int_{\{F_Q>\lambda\}}H_T\,d\lambda\\
 &\le [v]_{A_\infty(\D)} a^\eta v(Q)
   +\frac{2^n[v]_{A_\infty(\D)}\eta}{1+\eta}\int_Q H_T^{1+\eta}.
\end{align*}
We absorb the last term when $2^{n+2}[v]_{A_\infty(\D)}\eta\le1$. Thus $\int_Q H_T^{1+\eta}\le2[v]_{A_\infty(\D)} a^\eta v(Q)$. Applying the same layer-cake formula with measure $v(x)\,dx$ and using \eqref{eq:RHdist2},
\begin{align*}
 \int_Q v H_T^\eta
 &\le a^\eta v(Q)+\frac{2^n\eta}{1+\eta}\int_Q H_T^{1+\eta}\\
 &\le (1+2^{n+1}\eta [v]_{A_\infty(\D)})a^\eta v(Q)
 \le 2 a^\eta v(Q).
\end{align*}
Let $T\to\infty$. Since $v\le F_Q$ almost everywhere, monotone convergence gives the second inequality in \eqref{eq:dyadic-RH}. All cubes used in the proof belong to $\D$.
\end{proof}

\subsection{The matrix kernel and the local norm}\label{sub:dyadicRH}
\begin{proof}[Proof of \Lref{lem:kernel}]
Let $U=W^{-1/(p-1)}$ and $V'_Q=V_{U,p',Q}$. By \Lref{lem:domain}, $U$ is locally integrable for either family $\mathcal G$. For $(Z,s)=(W,p)$ or $(U,p')$ and any fixed $d\times d$ complex matrix $B$, one has
\begin{equation}\label{eq:columns}
 \av{\|Z^{1/s}B\|_{\op}}_{s,Q}\eqsim_{d,s}\|V_{Z,s,Q}B\|_{\op}.
\end{equation}
Indeed, any $d\times d$ matrix $C$ satisfies $\max_{1\le k\le d}|Ce_k|\le\|C\|_{\op}\le\sum_{k=1}^d |Ce_k|$. Apply this to $C=Z^{1/s}(y)B$. The $L^s$ average of each column length of $Z^{1/s}B$ is comparable to $|V_{Z,s,Q}Be_k|$ by \eqref{eq:reducing}. Summing over $k$ proves \eqref{eq:columns}.

By \eqref{eq:columns} and the equality of the operator norms of a matrix and its adjoint, the quantity defining the matrix characteristic on each $Q\in\mathcal G$ satisfies
\begin{align*}
 &\left(\fint_Q\left(\fint_Q\|W(x)^{1/p}W(y)^{-1/p}\|_{\op}^{p'}dy\right)^{p/p'}dx\right)^{1/p}\\
 &\eqsim_{d,p}\left(\fint_Q\|V'_QW(x)^{1/p}\|_{\op}^{p}dx\right)^{1/p}\\
 &=\av{\|W^{1/p}V'_Q\|_{\op}}_{p,Q}
 \eqsim_{d,p}\|V_QV'_Q\|_{\op}.
\end{align*}
The first comparison applies \eqref{eq:columns} with $(Z,s)=(U,p')$ to the inner integral, and the last applies it with $(Z,s)=(W,p)$ to the outer integral. The same calculation for $(U,p')$ gives the comparison with the dual characteristic, first allowing it to be infinite. The product of the reducing matrices is finite and uniformly bounded on $\mathcal G$, so $U\in A_{p'}(\mathcal G)$. Taking suprema yields
\begin{equation}\label{eq:product}
 [W]_{A_p(\mathcal G)}^{1/p}\eqsim_{d,p}\sup_{Q\in\mathcal G}\|V_QV'_Q\|_{\op}
 =\sup_{Q\in\mathcal G}\|V'_QV_Q\|_{\op}\eqsim_{d,p}[U]_{A_{p'}(\mathcal G)}^{1/p'}.
\end{equation}
In particular, for $K_Q=\|V_QW^{-1/p}\|_{\op}$,
\begin{equation}\label{eq:baseR}
 \av{K_Q}_{p',Q}=\av{\|U^{1/p'}V_Q\|_{\op}}_{p',Q}
 \lesssim_{d,p}\|V'_Q V_Q\|_{\op}\lesssim_{d,p}[W]_{A_p(\mathcal G)}^{1/p}.
\end{equation}

Let $e\ne0$ and $v_e=|U^{1/p'}e|^{p'}$. The inequality
\[
 \left(v_e(x)/v_e(y)\right)^{1/p'}
 \le\|U(x)^{1/p'}U(y)^{-1/p'}\|_{\op}
\]
shows, for every cube $Q$, that
\begin{align*}
 v_e(x)\left(\av{v_e^{-1/(p'-1)}}_Q\right)^{p'-1}
 &=\left(|U^{1/p'}(x)e|^p
       \av{|U^{1/p'}(y)e|^{-p}}_Q\right)^{1/(p-1)}\\
 &\le\left(\av{\|U(x)^{1/p'}U(y)^{-1/p'}\|_{\op}^{p}}_Q\right)^{p'/p}.
\end{align*}
Averaging in $x$ and taking the supremum over $Q\in\mathcal G$ gives $[v_e]_{A_{p'}(\mathcal G)}\le[U]_{A_{p'}(\mathcal G)}$. For the full-cube case, set $\mathcal G=\mathcal Q$. Applying the scalar Fujii--Wilson estimate and the sharp reverse H\"older theorem of Hyt\"onen and P\'erez~\cite[Proposition~2.2 and Theorem~2.3]{HP} yields the following bounds. For the first, their comparison with the exponential $A_\infty$ characteristic is followed by Jensen's inequality $\av{v}_Q\exp(\av{\log(1/v)}_Q)\le\av{v}_Q\av{v^{-1/(p'-1)}}_Q^{p'-1}$ at exponent $p'$.
\begin{align*}
 1\le [U]_{A_{p',\infty}^{\mathrm{sc}}}&\lesssim_{n,p}[U]_{A_{p'}}\lesssim_{n,d,p}[W]_{A_p}^{1/(p-1)},\\
 \av{v_e^{1+\delta}}_Q^{1/(1+\delta)}&\le2\av{v_e}_Q,
 \qquad\delta=c_n/[U]_{A_{p',\infty}^{\mathrm{sc}}}.
\end{align*}
Scalar multiples of $e$ have the same characteristic. In particular, the same $\delta$ works for the directions $V_Qe_k$, for every $Q$ and $k$. Taking adjoints of $V_QW^{-1/p}$, then applying the column inequality with $s=p'(1+\delta)$, gives
\begin{align*}
 \av{K_Q}_{s,Q}
 &\le\sum_{k=1}^d\av{|W^{-1/p} V_Qe_k|}_{s,Q}
 \le 2^{1/p'}\sum_{k=1}^d\av{|W^{-1/p}V_Qe_k|}_{p',Q}\\
 &\le 2^{1/p'}d\,\av{K_Q}_{p',Q}.
\end{align*}
Equation \eqref{eq:baseR} concludes the full-cube estimate.

For the dyadic assertions, restrict every cube in \eqref{eq:columns}--\eqref{eq:baseR} to $\D$. The argument bounding $[v_e]_{A_{p'}}$ now gives $[v_e]_{A_{p'}(\D)}\le[U]_{A_{p'}(\D)}$. By \Lref{lem:dyadicRH},
\[
 [U]_{A_{p',\infty}^{\mathrm{sc}}(\D)}\lesssim_p[U]_{A_{p'}(\D)}\lesssim_{d,p}[W]_{A_p(\D)}^{1/(p-1)}.
\]
Use the second part of \eqref{eq:dyadic-RH} with $\eta=c_n/[U]_{A_{p',\infty}^{\mathrm{sc}}(\D)}$, choosing $c_n$ sufficiently small for both the full and the dyadic estimates. Each column satisfies the same reverse H\"older estimate as above. The column sum then proves \eqref{eq:kernel} with the restricted characteristics. This proves the restricted version of \Lref{lem:kernel}.
\end{proof}

\begin{proof}[Proof of \Lref{lem:localnorm}]
Fix $Q\in\mathcal G$ and $W\in A_p(\mathcal G)$, with $\mathcal G=\mathcal Q$ or a single grid. The kernel $K_Q$ belongs to $L^{p'}(Q)$ by \eqref{eq:domainlocal} and the local integrability of $W^{-1/p}$. For every sequence of length $N\in\mathbb N$, write $F(y)=\|(h_j(y))_j\|_{\ell_{N,d}^q}$. Minkowski's inequality for the sequence norm and H\"older's inequality give
\begin{align*}
 \|\overrightarrow B_Q\mathbf h\|_{L^p(\ell^q)}
 &= |Q|^{1/p}\left\|(\av{|V_Q W^{-1/p}h_j|}_Q)_j\right\|_{\ell^q}\\
 &\le |Q|^{1/p}\fint_Q K_Q(y)F(y)\,dy\\
 &\le \av{K_Q}_{p',Q}\|\mathbf h\|_{L^p(\ell_d^q)}.
\end{align*}
For $q=\infty$ the same bound uses the supremum in place of the sum. Countable sequences follow by truncation. For the reverse bound choose a countable dense set $(u_k)_{k\ge1}$ on the Euclidean unit sphere. Given $0<\eta<1$, choose $u_\eta(y)$ to be the first $u_k$ satisfying
$|V_QW(y)^{-1/p}u_k|>(1-\eta)K_Q(y)$. This is a measurable choice because the functions in the inequality are measurable. The choice is possible wherever $W(y)$ is positive definite, by continuity in the unit vector. Define a single function on $Q$, extended by zero off $Q$, by
\[
 h_\eta(y)=K_Q(y)^{p'-1}u_\eta(y).
\]
Since $(p'-1)p=p'$, we obtain
\begin{align*}
 \|h_\eta\|_{L^p(\mathbb C^d)}^p&=\int_Q K_Q(y)^{p'}\,dy,\\
 \|B_Q h_\eta\|_p&\ge(1-\eta)|Q|^{1/p}\fint_Q K_Q(y)^{p'}\,dy.
\end{align*}
Consequently $\|B_Q\|_{L^p(\mathbb C^d)\to L^p}\ge(1-\eta)\av{K_Q}_{p',Q}$. Letting $\eta\downarrow0$ and taking one nonzero entry in the sequence proves \eqref{eq:local-story}.

For $W=wI_d$ with $V_Q=\av{w}_Q^{1/p}I_d$ and $\sigma=w^{-1/(p-1)}$, we have $K_Q=(\av{w}_Q/w)^{1/p}$. Substitution into \eqref{eq:local-story} gives \eqref{eq:one-cube-cost}, with equality for $h=\sigma^{1/p} e\one_Q$.

Combining the exact identity \eqref{eq:local-story} with the column comparison \eqref{eq:columns}, we obtain for every $Q\in\mathcal G$
\[
 \av{K_Q}_{p',Q}\eqsim_{d,p}\|V'_QV_Q\|_{\op}.
\]
Taking the supremum over $Q\in\mathcal G$ and using \eqref{eq:product} yields
\[
 \sup_{Q\in\mathcal G}\|\overrightarrow B_Q\|_{L^p(\ell_d^q)\to L^p(\ell^q)}\eqsim_{d,p}[W]_{A_p(\mathcal G)}^{1/p}.
\]
\end{proof}

\subsection{The common sparse family}\label{subsec:selection-proof}
\begin{proof}[Proof of \Lref{lem:weak}]
For $q=\infty$, use $\sup_jM_Q^{\mathcal U}u_j\le M_Q^{\mathcal U}F$ and the scalar weak $(1,1)$ estimate. Let $q<\infty$ and set $a=\lambda/(2^nq')$. If $\av F_Q>a$, the claim follows from $|Q|<a^{-1}\int_QF$. Otherwise, let $\mathcal B$ be the maximal cubes of $\mathcal U$ with $\av F_B>a$, and let $E=\bigcup_{B\in\mathcal B}B$. Then
\[
 |E|\le a^{-1}\int_Q F,\qquad
 a<\av F_B\le 2^na=\lambda/q'.
\]
Let $\mathcal P$ consist of the cubes in $\mathcal B$ and all terminal cubes of $\mathcal U$ not contained in $E$. This is a partition of $Q$. Define
\[
 h_j=\sum_{P\in\mathcal P}\av{u_j}_P\one_P,\qquad H=\|(h_j)_j\|_{\ell^q}.
\]
Every $P\in\mathcal P$ satisfies $\av F_P\le\lambda/q'$, by the preceding bound on stopped cubes and by $\av F_P\le a$ on the remaining terminal cubes. Hence
\[
 \|H\|_{L^\infty(Q)}\le\lambda/q',\qquad \int_Q H\le\int_Q F.
\]
For $x\in Q\setminus E$ and $R\in\mathcal U$ containing $x$, the dyadic nesting property gives $R=\bigcup_{P\in\mathcal P,\,P\subseteq R}P$. Indeed, a stopped cube meeting $R$ is contained in $R$, and a terminal cube cannot properly contain $R$. Thus
\[
 \int_R h_j=\int_R u_j,\qquad M_Q^{\mathcal U}h_j(x)=M_Q^{\mathcal U}u_j(x).
\]
The scalar dyadic maximal inequality $\|M^\D h\|_q\le q'\|h\|_q$ gives
\begin{align*}
 |\{x\in Q\setminus E:\|(M_Q^{\mathcal U}u_j(x))_j\|_{\ell^q}>\lambda\}|
 &\le\lambda^{-q}\sum_j\|M_Q^{\mathcal U}h_j\|_{L^q(Q)}^q\\
 &\le\lambda^{-q}(q')^q\int_Q H^q\\
 &\le\frac{q'}{\lambda}\int_Q F.
\end{align*}
Adding $|E|\le 2^nq'\lambda^{-1}\int_Q F$ proves the lemma.
\end{proof}

\begin{proof}[Proof of \Pref{prop:freeze}]
For $Q\in\mathcal T$, write
\begin{equation}\label{eq:freeze-data}
 \Phi_Q=\|W^{1/p}V_Q^{-1}\|_{\op}^p,\qquad
 u_j=|V_Qf_j|,\qquad m_Q=\av{\|(u_j)_j\|_{\ell^q}}_Q.
\end{equation}
The lower reducing inequality on $Q$ and the upper one on $R\subseteq Q$ give
\begin{align*}
 \left(\fint_Q|W^{1/p}V_Q^{-1}e_k|^p\right)^{1/p}&\le1,\quad 1\le k\le d,\\
 |V_RV_Q^{-1}v|&\le c_d\left(\fint_R|W^{1/p}V_Q^{-1}v|^p\right)^{1/p}
 \quad (|v|=1).
\end{align*}
Summing the column bounds in the first line and taking the supremum in the second gives
\begin{equation}\label{eq:relfreeze}
 \av{\Phi_Q}_Q^{1/p}\le d,\qquad
 \|V_RV_Q^{-1}\|_{\op}\le c_d\av{\Phi_Q}_R^{1/p}.
\end{equation}
The second estimate is uniform over all unit vectors.

Set $L=4d^p$ and $b_q=4(2^n+1)q'$. If $m_Q=0$, all $f_j$ vanish almost everywhere on $Q$ and we take $\mathcal C(Q)=\varnothing$. For $m_Q>0$, select the maximal strict descendants $P$ of $Q$ in $\mathcal T$ satisfying at least one of
\begin{equation}\label{eq:badset}
 \av{\Phi_Q}_P>L,
 \qquad
 \left\|\left(\max_{\substack{R\in\mathcal T\\P\subseteq R\subseteq Q}}\av{u_j}_R\right)_j\right\|_{\ell^q}>b_qm_Q.
\end{equation}
For $m_Q>0$, neither condition holds at $Q$. The selected descendants, denoted by $\mathcal C(Q)$, are disjoint. On their union, either $M_Q^{\mathcal T}\Phi_Q>L$ or $\|(M_Q^{\mathcal T}u_j)_j\|_{\ell^q}>b_qm_Q$. By \eqref{eq:relfreeze}, the scalar weak $(1,1)$ bound and \Lref{lem:weak},
\begin{align*}
 \sum_{P\in\mathcal C(Q)}|P|
 &\le \frac{1}{L}\int_Q\Phi_Q
      +\frac{(2^n+1)q'}{b_qm_Q}\int_Q\|(u_j)_j\|_{\ell^q}\\
 &\le\frac{|Q|}{4}+\frac{|Q|}{4}=\frac{|Q|}{2}.
\end{align*}
For $m_Q=0$ the sum is zero. Iterating from the roots gives $\Spar$ and the sets
\[
 E_Q=Q\setminus\bigcup_{P\in\mathcal C(Q)}P,\qquad Q\in\Spar,
\]
are disjoint and have measure at least $|Q|/2$.

For each $Q\in\Spar$, consider the cubes assigned to $Q$, namely those for which $Q$ is the smallest selected ancestor:
\[
 \mathcal T_Q=\{R\in\mathcal T:R\subseteq Q,\ R\not\subseteq P\text{ for every }P\in\mathcal C(Q)\}.
\]
The groups partition $\mathcal T$. Fix $x\in Q$ and suppose $m_Q>0$. The cubes in $\mathcal T_Q$ containing $x$ form the chain from $Q$ to a smallest cube $S$. Every cube $R$ in this chain fails the first stopping condition, while $S$ fails the second. Thus
\begin{align*}
 \av{|V_Rf_j|}_R
 &\le\|V_RV_Q^{-1}\|_{\op}\av{|V_Qf_j|}_R
 \le c_d L^{1/p}\av{u_j}_R,\\
 \left\|\left(\max_{\substack{R\in\mathcal T\\S\subseteq R\subseteq Q}}\av{u_j}_R\right)_j\right\|_{\ell^q}
 &\le b_qm_Q.
\end{align*}
Combining the two inequalities gives
\begin{equation}\label{eq:commonpoint}
 \left\|\left(\max_{\substack{R\in\mathcal T_Q\\x\in R}}\av{|V_R f_j|}_R\right)_j\right\|_{\ell^q}
 \le c_dL^{1/p} b_qm_Q=C_{n,d,p}q'm_Q.
\end{equation}
If $m_Q=0$ all descendant averages vanish, and if $x\notin Q$ the maxima are zero.

Because the groups partition the forest, for $1<q<\infty$ we obtain
\begin{align}\label{eq:recursion}
 \sum_{j=1}^N\left(\max_{\substack{R\in\mathcal T\\x\in R}}\av{|V_R f_j|}_R\right)^q
 &\le\sum_{Q\in\Spar}\sum_{j=1}^N\left(\max_{\substack{R\in\mathcal T_Q\\x\in R}}\av{|V_R f_j|}_R\right)^q\notag\\
 &\le (C_{n,d,p}q')^q\sum_{Q\in\Spar}\one_Q(x)m_Q^q.
\end{align}
For $q=\infty$, the maxima over $j$ and $Q\in\Spar$ commute, giving
\[
 \max_{1\le j\le N}\max_{\substack{R\in\mathcal T\\x\in R}}\av{|V_R f_j|}_R
 \le C_{n,d,p}\max_{Q\in\Spar}\one_Q(x)m_Q.
\]
Taking $q$th roots in \eqref{eq:recursion}, or the supremum when $q=\infty$, and substituting $m_Q$ from \eqref{eq:freeze-data} proves \eqref{eq:sparsefreeze}.
\end{proof}

\subsection{The scalar sparse sum and the power weights}\label{subsec:scalar-proof}
\begin{proof}[Proof of \Lref{lem:scalar-sparse}]
Since the entries are nonnegative, it suffices to estimate $S=\sum_{Q\in\Spar}\one_Q\av F_{r,Q}$ for every $q$. Let $g\ge0$, $\|g\|_{p'}=1$. On the disjoint sparse sets $E_Q$, we have $\av F_{r,Q}\le H$ and $\av g_Q\le M^\D g$. Hence
\begin{align*}
 \int_{\R^n}Sg
 &=\sum_{Q\in\Spar}|Q|\av F_{r,Q}\av g_Q
 \le 2\sum_{Q\in\Spar}\int_{E_Q}H M^\D g\\
 &\le2\|H\|_p\|M^\D g\|_{p'}\le2p\|H\|_p.
\end{align*}
Duality gives \eqref{eq:scalar-sparse} with $C_p=2p$ for every $q$.
\end{proof}

\begin{proof}[Proof of \Lref{lem:power}]
On $I=(0,1)$, direct integration gives
\[
 \av{|x|^\alpha}_I=\frac1{1+\alpha},\qquad
 \av{|x|^{-\alpha/(r-1)}}_I=\frac1{1-\alpha/(r-1)}.
\]
The product $\av{|x|^\alpha}_I\av{|x|^{-\alpha/(r-1)}}_I^{r-1}$ is at least one by H\"older's inequality. For an interval $J$, set $a=\operatorname{dist}(J,0)$. If $a\ge |J|$, then $a\le |x|\le 2a$ for $x\in J$. Hence
\[
 \av{|x|^\alpha}_J\av{|x|^{-\alpha/(r-1)}}_J^{r-1}
 \le 2^{2|\alpha|}\le 2^{2\max\{1,r-1\}}.
\]
If $a<|J|$, then $J\subseteq(-2|J|,2|J|)$. For $\beta>-1$,
\[
 \av{|x|^\beta}_J\le\frac{2^{\beta+2}}{1+\beta}|J|^\beta.
\]
Using $\beta=\alpha$ and $\beta=-\alpha/(r-1)$ gives
\[
 \av{|x|^\alpha}_J\av{|x|^{-\alpha/(r-1)}}_J^{r-1}
 \le \frac{2^{2r}}{1+\alpha}\left(1-\frac{\alpha}{r-1}\right)^{1-r}.
\]
The corresponding cube averages reduce to interval averages in the first coordinate. The lower test on $(0,1)^n$ and the two upper estimates prove \eqref{eq:powerformula}.
\end{proof}

\Needspace{5\baselineskip}
\smallskip
\noindent{\bf Funding }
Not applicable.
\smallskip

\noindent{\bf Data Availability} This manuscript has no associated data.

\smallskip
\noindent{\bf\Large Declarations}
\smallskip

\noindent\textbf{Conflicts of Interest}\space The authors declare that there are no conflicts of interest.

\smallskip
\noindent{\bf AI disclosure.}
During the preparation of this manuscript, the authors used ChatGPT and Codex to discuss possible proof strategies, assist in drafting and checking parts of the mathematical arguments, organize references, and improve the clarity and readability of the exposition. All mathematical statements and proofs in the final manuscript were written, reviewed, and verified by the authors, who take full responsibility for the content of this manuscript.

\end{document}